\documentclass[reqno,11pt]{article}
\usepackage{psfig, amsmath, amstexnb, amsthm}
\usepackage{amssymb}  
\newtheorem{theorem}{Theorem}[section]
\newtheorem{prop}[theorem]{Proposition}
\newtheorem{lemma}[theorem]{Lemma}
\newtheorem{cor}[theorem]{Corollary}
\newtheorem{definition}[theorem]{Definition}

\theoremstyle{definition}

\newtheorem{remark}[theorem]{Remark}
\newtheorem{assump}[theorem]{Assumption}

\numberwithin{equation}{section}

\makeatletter

\def\enum{\ifnum \@enumdepth >3 \@toodeep\else
        \advance\@enumdepth \@ne 
        \edef\@enumctr{enum\romannumeral\the\@enumdepth}\list
        {\csname label\@enumctr\endcsname}
        {\setlength{\topsep}{1mm}
        \setlength{\parsep}{0mm}
        \setlength{\itemsep}{0mm}
        \setlength{\labelsep}{2mm}
        \settowidth{\leftmargin}{M.}
        \addtolength{\leftmargin}{\labelsep}
        \usecounter{\@enumctr}
        \def\makelabel##1{\hss\llap{##1}}}\fi}

\def\itemiz{\ifnum \@itemdepth >3 \@toodeep\else \advance\@itemdepth \@ne
        \edef\@itemitem{labelitem\romannumeral\the\@itemdepth}%
        \list{\csname\@itemitem\endcsname}{
        \setlength{\topsep}{1mm}
        \setlength{\parsep}{0mm}
        \setlength{\parsep}{0mm}
        \setlength{\itemsep}{0mm}
        \setlength{\labelsep}{2mm}
        \settowidth{\leftmargin}{M.}
        \addtolength{\leftmargin}{\labelsep}
        \def\makelabel##1{\hss\llap{##1}}}\fi}

\def\captionheadfont@{\scshape}
\def\captionfont@{\small}
\long\def\@makecaption#1#2{%
  \setbox\@tempboxa\vbox{\color@setgroup
    \advance\hsize-3pc\noindent
    \captionfont@\captionheadfont@#1\@xp\@ifnotempty\@xp
        {\@cdr#2\@nil}{.\captionfont@\upshape\enspace#2}%
    \unskip\kern-3pc\par
    \global\setbox\@ne\lastbox\color@endgroup}%
  \ifhbox\@ne 
    \setbox\@ne\hbox{\unhbox\@ne\unskip\unskip\unpenalty\unkern}%
  \fi
  \ifdim\wd\@tempboxa=\z@ 
    \setbox\@ne\hbox to\columnwidth{\hss\kern-3pc\box\@ne\hss}%
  \else 
    \setbox\@ne\vbox{\unvbox\@tempboxa\parskip\z@skip
        \noindent\unhbox\@ne\advance\hsize-3pc\par}%
\fi
  \ifnum\@tempcnta<64 
    \addvspace\abovecaptionskip
    \moveright 1.5pc\box\@ne
  \else 
    \moveright 1.5pc\box\@ne
    \nobreak
    \vskip\belowcaptionskip
  \fi
\relax
}

\def\makeoverbar#1#2#3#4#5#6#7#8{{%
 \setbox0=\hbox{$\m@th#2\mkern#5mu{#3}\mkern#6mu$}%
 \setbox1=\null \dimen@=#4\fontdimen8#13 \dimen@=3\dimen@ 
 \advance\dimen@ by \ht0 \dimen@=-#7\dimen@ \advance\dimen@ by \wd0
 \wd1=\dimen@ \dp1=\dp0 
 \dimen@=#4\fontdimen8#13
 \dimen@i=\fontdimen8#13
 \fontdimen8#13=#8\dimen@
 \advance\dimen@ by -\fontdimen8#13 \dimen@=3\dimen@
 \advance\dimen@ by \ht0 \ht1=\dimen@ 
 \rlap{\hbox to \wd0{$\m@th\hss#2{\overline{\box1}}\mkern#5mu$}}
 \fontdimen8#13=\dimen@i}}

\makeatother

\DeclareMathSymbol{\leqsymb}{\mathalpha}{AMSa}{"36}
\def\leqs{\mathrel\leqsymb}
\DeclareMathSymbol{\geqsymb}{\mathalpha}{AMSa}{"3E}
\def\geqs{\mathrel\geqsymb}
\DeclareMathSymbol{\gtreqqlesssymb}{\mathalpha}{AMSa}{"54}

\def\,{\ifmmode\mskip\thinmuskip\else\kern\fontdimen3}

\newcommand{\field}[1]{\mathbb{#1}}

\newcommand{\N}{\field{N}\,}    
\newcommand{\E}{\field{E}}      
\newcommand{\fP}{\field{P}}     

\newcommand{\cA}{{\mathcal A}}  
\newcommand{\cF}{{\mathcal F}}  

\DeclareMathOperator{\e}{e}             
\DeclareMathOperator{\dd}{d}            
\DeclareMathOperator{\variance}{Var}    

\def\math#1{\ifmmode
\mathchoice{\mbox{$\displaystyle\rm#1$}}
{\mbox{$\textstyle\rm#1$}}
{\mbox{$\scriptstyle\rm#1$}}
{\mbox{$\scriptscriptstyle\rm#1$}}\else
{\mbox{$\rm#1$}}\fi}            
\DeclareMathOperator{\defby}{\raisebox{0.35pt}{\math{:}}\!\!=}

\def\defwd#1{{\em #1}}                  

\def\nbh{neighbourhood}

\def\eps{\varepsilon}

\def\w{\omega}

\def\z{\zeta}

\def\xdet{x^{\det}}

\def\xc{x_{\math c}}
\def\tc{t_{\math c}}
\def\tauc{\tau_{\math c}}
\def\lc{\lambda_{\math c}}
\def\xper{x^{\math{per}}}
\def\xperof#1{x^{{\math{per}},#1}}
\def\xdetof#1{x^{\det,#1}}
\def\esa0{(\eps\sqrt{a_0}\mskip1.5mu)}

\def\6#1{\dd\!#1}                       
\def\dpar#1#2{\frac{\partial #1}{\partial #2}}  
\def\dtot#1#2{\frac{\6{#1}}{\6{#2}}}  

\def\Bigevalat#1{\Bigr|_{#1}^{\phantom{#1}}}    

\def\bigpar#1{\bigl(#1\bigr)}                   
\def\biggpar#1{\biggl(#1\biggr)}        
\def\Bigpar#1{\Bigl(#1\Bigr)}

\def\bigbrak#1{\bigl[#1\bigr]}          
\def\Bigbrak#1{\Bigl[#1\Bigr]}          
\def\biggbrak#1{\biggl[#1\biggr]}

\def\set#1{\{#1\}}                              
\def\bigset#1{\bigl\{#1\bigr\}}         
\def\Bigset#1{\Bigl\{#1\Bigr\}}         
\def\biggset#1{\biggl\{#1\biggr\}}

\def\setsuch#1#2{\{#1\colon #2\}}                
\def\bigsetsuch#1#2{\bigl\{#1\colon #2\bigr\}}

\def\abs#1{\lvert#1\rvert}                      
\def\bigabs#1{\bigl|#1\bigr|}           
           
\def\biggabs#1{\biggl|#1\biggr|}

\def\Order#1{{\mathcal O}(#1)}                  
\def\bigOrder#1{{\mathcal O}\bigl(#1\bigr)}     
\def\BigOrder#1{{\mathcal O}\Bigl(#1\Bigr)}     
\def\biggOrder#1{{\mathcal O}\biggl(#1\biggr)}  

\def\bigprob#1{\fP\bigl\{#1\bigr\}}

\def\bigexpec#1{\E\bigl\{#1\bigr\}}

\def\indexfct#1{1_{\set{#1}}}           
\def\indicator#1{1_{#1}}                
\def\bigindexfct#1{1_{\bigset{#1}}}           

\def\bigprobin#1#2{\fP^{\mskip1.5mu #1}\bigl\{#2\bigr\}}    
\def\Bigprobin#1#2{\fP^{\mskip1.5mu #1}\Bigl\{#2\Bigr\}} 
\def\biggprobin#1#2{\fP^{\mskip1.5mu #1}\biggl\{#2\biggr\}} 

\def\bigexpecin#1#2{\E^{\mskip1.5mu #1}\bigl\{#2\bigr\}}    
\def\Bigexpecin#1#2{\E^{\mskip1.5mu #1}\Bigl\{#2\Bigr\}}    
\def\biggexpecin#1#2{\E^{\mskip1.5mu #1}\biggl\{#2\biggr\}} 

\def\bibtitle#1#2{#1, {\em #2}}                       
\def\bibref#1#2#3#4#5{#1 {\bf #2}:#3--#4 (#5)}        
\def\bibarticle#1#2#3#4#5#6#7{\bibtitle{#1}{#2},
\bibref{#3}{#4}{#5}{#6}{#7}.}

\def\bibbook#1#2#3#4{#1, {\em #2} (#3, #4).}

\def\CMP{Comm.\ Math.\ Phys.}

\def\JPA{J.\ Phys.\ A}

\def\JSP{J.\ Stat.\ Phys.}
\def\Nat{Nature}

\def\PRA{Phys.\ Rev.\ A}
\def\PRB{Phys.\ Rev.\ B}
\def\PRE{Phys.\ Rev.\ E}
\def\PRL{Phys.\ Rev.\ Letters}

\begin{document}


\title{The effect of additive noise on dynamical hysteresis}
\author{Nils Berglund and Barbara Gentz}
\date{}

\maketitle

\begin{abstract}
\noindent
We investigate the properties of hysteresis cycles produced by a
one-dimensional, periodically forced Langevin equation. We show that
depending on amplitude and frequency of the forcing and on noise intensity,
there are three qualitatively different types of hysteresis cycles. Below a
critical noise intensity, the random area enclosed by hysteresis cycles is
concentrated near the deterministic area, which is different for small and
large driving amplitude. Above this threshold, the area of typical
hysteresis cycles depends, to leading order, only on the noise intensity.
In all three regimes, we derive mathematically rigorous estimates for
expectation, variance, and the probability of deviations of the
hysteresis area from its typical value.  
\end{abstract}

\leftline{\small{\it Date.\/} July 27, 2001.}
\leftline{\small 2000 {\it MSC.\/} 
37H20 (primary), 60H10, 34C55, 34E15, 82C31 (secondary).}


\section{Introduction}
\label{sec_in}

For a long time, hysteresis was considered as a purely static phenomenon.
As a consequence, it has been modeled by various integral operators
relating the \lq\lq output\rq\rq\ of the system to its \lq\lq input\rq\rq,
for operators not depending on the speed of variation of the input (see for
instance \cite{May} and \cite{MNZ} for reviews). 

This situation changed drastically a decade ago, when Rao and coauthors
published a numerical study of the effect of the input's frequency on shape
and area of hysteresis cycles \cite{RKP}. They proposed in particular that
the area $\cA$ of a hysteresis cycle, which measures the energy dissipation
per period, should obey a scaling law of the form 
\begin{equation}
\label{in1}
\cA \simeq A^\alpha \eps^\beta
\end{equation}
for small amplitude $A$ and frequency $\eps$ of the periodic input
(e.\,g.\ the magnetic field), and some model-dependent exponents $\alpha$ and
$\beta$. This work triggered a substantial amount of numerical, experimental
and theoretical studies, trying to establish the validity of the scaling law
\eqref{in1} for various systems, a problem which has become known as {\em
dynamical hysteresis}. 

The first model investigated in \cite{RKP} is a Langevin partial
differential equation for the spatially extended, $N$-component order
parameter (e.\,g.\ the magnetization), in a $(\Phi^2)^2$-potential
with $O(N)$-symmetry, in the limit $N\to\infty$. Their numerical
experiments suggested that \eqref{in1} holds with $\alpha\simeq2/3$
and $\beta\simeq1/3$. Various theoretical arguments \cite{DT,SD,ZZ}
indicate that the scaling law should be valid, but with $\alpha=\beta=1/2$. 

The second model considered in \cite{RKP} is an Ising model with
Monte-Carlo dynamics. Here the situation is not so clear.  Different
numerical  simulations (for instance \cite{LP,AC,ZZL}) suggested scaling
laws with widely different exponents. More careful simulations \cite{SRN},
however,  showed that the behaviour of hysteresis cycles depends in a
complicated way on the mechanism of magnetization reversal, and no
universal scaling law of the form \eqref{in1} should be expected.  Rigorous
results on hysteresis in the Ising model are only available for
discontinuous reversal (quenching) of the field \cite{SS}.

A third kind of models for which scaling laws of hysteresis cycles have been
investigated belong to the mean field class, and include the
Curie$\mskip1.5mu$--$\mskip-1.5mu$Weiss 
model. A one-dimensional deterministic equation modeling a bistable laser,
and being equivalent to the equation of motion of an overdamped particle in a
periodically forced double-well potential,  was considered in
\cite{Jung}. The area of hysteresis cycles was shown to obey the scaling law
\begin{equation}
\label{in2}
\cA \simeq \cA_0 + \eps^{2/3}
\end{equation}
for sufficiently large driving amplitude. A similar equation governs the
dynamics of the magnetization in the
Curie$\mskip1.5mu$--$\mskip-1.5mu$Weiss model, in the limit of 
infinite system size. This equation was examined in \cite{TO}, where it was
shown that the behaviour changes drastically when the amplitude of the
forcing crosses a threshold, a phenomenon they termed \lq\lq dynamic phase
transition\rq\rq. 

As pointed out in \cite{Ra}, the difference between the scaling laws
\eqref{in1} and \eqref{in2} can be attributed to the existence of a
potential barrier for the one-dimensional order parameter, which is absent
in higher dimensions. The deterministic equation, however, neglects both
thermal fluctuations and the finite system size, whose effects may be
modeled by an additive white noise (see for instance \cite{Martin}). Noise,
however, may help to overcome the potential barrier and change the scaling
law. 

The aim of the present work is to give a rigorous characterization of the
effect of additive white noise on scaling properties of hysteresis cycles.
For definiteness, we shall consider the case of a Ginzburg--Landau
potential, i.\,e., the stochastic differential equation 
\begin{equation}
\label{in3}
\6x_s = -\dpar{}x \Bigbrak{\frac14 x_s^4 - \frac12 x_s^2 - \lambda(\eps
s)x_s} \6s + \sigma \6W_s,
\end{equation}
where $W_s$ is a standard Brownian motion, and 
\begin{equation}
\label{in4}
\lambda(\eps s) = -A\cos(2\pi\eps s), 
\quad A>0. 
\end{equation}
However, our results depend only on certain qualitative features of the
bifurcation diagram and the proofs carry over to a more general setup as
in~\cite{BG1,BG2}.

In the deterministic case $\sigma=0$, it is known \cite{TO,Jung,BK} that 
\begin{itemiz}
\item	for $A<\lc+\Order{\eps}$ (where $\lc=2/(3\sqrt3)$ is such that the
potential has two wells if and only if $\abs{\lambda}<\lc$), solutions of
\eqref{in3} are attracted by hysteresis cycles (one for each potential well)
enclosing an area of order $\eps$, and with nonzero mean;

\item	for $A>\lc+\Order{\eps}$, solutions are attracted by a hysteresis
cycle enclosing an area of order $\cA_0+\eps^{2/3}(A-\lc)^{1/3}$, where the
static hysteresis area $\cA_0$ is a constant, depending only on the geometry
of the equilibrium branches. 
\end{itemiz}

\begin{figure}
 \centerline{\psfig{figure=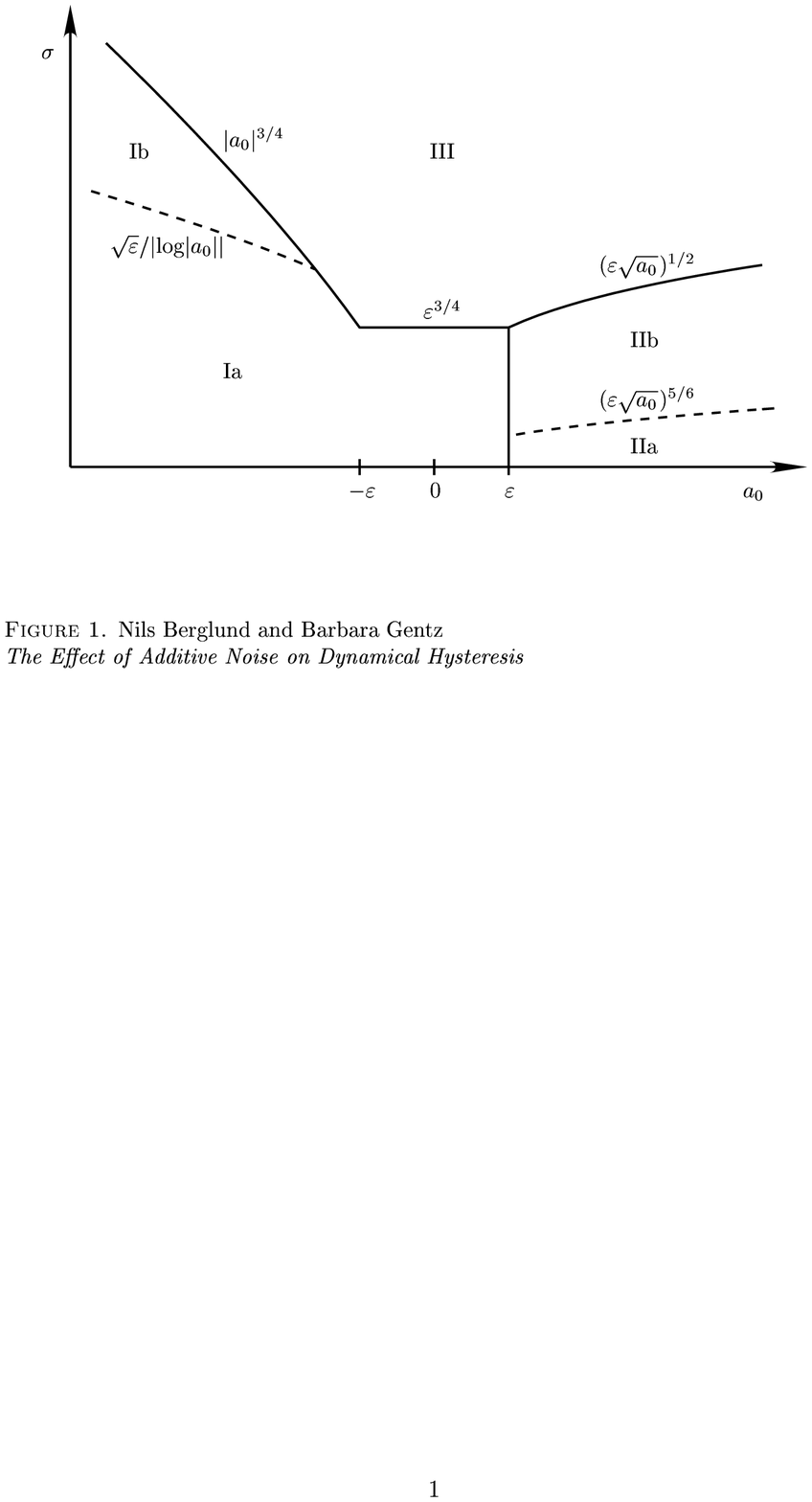,width=120mm,clip=t}}
 \caption[]
 {The different regimes as a function of amplitude $A=\lc+a_0$ and noise
 intensity $\sigma$, for a given value of the frequency $\eps$.}
\label{fig1}
\end{figure}

For positive $\sigma$, the area $\cA$ enclosed by a trajectory during one
period is a random variable, depending on the realization $W_s(\w)$ of the
Brownian motion. Our aim is to characterize the distribution of $\cA$ as a
function of the parameters $\eps$, $\sigma$ and $a_0 = A-\lc$. It turns out
that the distribution is usually concentrated around a deterministic
reference value. We determine the expectation and variance of $\cA$.
Furthermore, we estimate the behaviour of deviations of $\cA$ from its
reference value. 

One of the main results is the existence of a threshold value for the noise
intensity $\sigma$, depending on $A$ and $\eps$: Below this threshold, the
area is concentrated near the corresponding deterministic value, while above
the threshold, it depends, to leading order, only on the noise intensity and
is slightly smaller than $\cA_0$. 

There are thus three parameter regimes, as shown in Figure~\ref{fig1}, with
qualitatively different behaviour of the area $\cA$. 

\begin{itemiz}
\item	In Case I, the {\em small amplitude regime}, the area is close to the
deterministic value of order $\eps$. There is a further subdivision into
Case Ia, where the distribution of $\cA$ is close to a Gaussian with
standard deviation $\sigma\sqrt{\eps}$ smaller than $\eps$, and Case Ib,
where the distribution is more spread out (see
Theorem~\ref{thm_smallamp} and Figure~\ref{fig3}). 

\item	In Case II, the {\em large amplitude regime}, the area is concentrated
near the deterministic value of order $\cA_0 + \esa0^{2/3}$. In Case
IIa, the distribution is close to a Gaussian with standard deviation of
order $\sigma \esa0^{1/6}$. In Case IIb, we can only show that $\cA$ is
concentrated in an interval of width $\esa0^{2/3}$ (see
Theorem~\ref{thm_largeamp} and Figure~\ref{fig4}).

\item	In Case III, the {\em  large noise regime}, $\cA$ is likely to
be close to a reference area $\hat\cA$ of order $\cA_0-\sigma^{4/3}$,
which is smaller 
than the static hysteresis area. This is due to the noise driving $x$ over
the potential barrier before it becomes minimal or vanishes. The deviation
$-\sigma^{4/3}$ does not depend on $\eps$ or $A$ (see
Theorem~\ref{thm_largenoise} and Figure~\ref{fig5}). 
\end{itemiz}

Hysteresis does not only occur in ferromagnets and lasers, but also in
mechanical systems displaying relaxation oscillations, such as the Van der
Pol oscillator. Here additive noise can also have the
effect of enabling jumps between stable states separated by a potential
barrier~\cite{Freidlin}. Simple climate models can also display
hysteresis, as has been observed for instance for the Atlantic
thermohaline circulation \cite{Rah,Monahan}. In 
these systems, the effect of small scale degrees of freedom is represented
by additive noise. Our results describe quantitatively how noise may cause
the system to switch to another equilibrium state, at an earlier time than
expected from the deterministic approximation. 

We presented our results in detail in Section~\ref{sec_res}.
Section~\ref{sec_det} contains a short description of the deterministic
dynamics, while the remaining sections present the proofs for the various
parameter regimes. 

\subsubsection*{Acknowledgements:}
B.\,G. thanks the Forschungsinstitut f\"ur Mathematik at ETH Z\"urich
and its director Professor Marc Burger for kind hospitality.


\section{Results}
\label{sec_res}

We consider the non-autonomous SDE 
\begin{equation}
\label{SDE0}
\6x_s = F(x_s,\lambda(\eps s)) \6s + \sigma \6W_s,
\end{equation}
where $F$ derives from a periodically forced double-well potential and $W_s$
is a standard Brownian motion on some probability space $(\Omega, \cF, \fP)$.
For definiteness, we shall consider the case
\begin{align}
\label{r1a}
F(x,\lambda) & = x - x^3 + \lambda 
= -\dpar{}x \Bigbrak{\frac14 x^4 - \frac12 x^2 - \lambda x} \\ 
\label{r1b}
\lambda(\eps s) & = -A\cos(2\pi\eps s), 
\quad A>0.
\end{align}
We introduce the notation $\fP^{\mskip1.5mu t_0,x_0}$ for the law of
the process $\{x_t\}_{t\geqs t_0}$, starting in $x_0$ at time $t_0$,
and use $\E^{\mskip1.5mu t_0,x_0}$ to denote expectations with respect
to $\fP^{\mskip1.5mu t_0,x_0}$. Note that the stochastic process
$\{x_t\}_{t\geqs t_0}$ is an inhomogeneous Markov process. 

Before turning to the precise statements of our results, let us
introduce some notations. We shall use
\begin{itemiz}
\item 
$y\vee z$ and $y \wedge z$ to denote the maximum or minimum,
respectively, of two real numbers $y$ and $z$.
\item
If $\varphi(t,\eps)$ and $\psi(t,\eps)$ are defined for small $\eps$ and
for $t$ in a given interval $I$, we write
$\psi(t,\eps)\asymp\varphi(t,\eps)$ if there exist strictly positive
constants $c_\pm$ such that $c_-\varphi(t,\eps) \leqs \psi(t,\eps) \leqs
c_+\varphi(t,\eps)$ for all $t\in I$ and all sufficiently small
$\eps$. The constants $c_\pm$ are understood to be independent of $t$
and $\eps$ (and hence also independent of small quantities like $\sigma$
and, possibly, $a_0$, which we consider as functions of $\eps$). 
\item 
By $g(u)=\Order{u}$ we indicate that there exist $\delta>0$ and $K>0$
such that $g(u)\leqs K u$ for all $u\in[0,\delta]$, where $\delta$ and
$K$ of course do not depend on $\eps$ or on the other small parameters
$a_0$ and $\sigma$. 
\item
Let $I$ be an interval.  The notation $\indicator{I}(x)$ is used for the
indicator function, taking value $1$ if $x\in I$ and $0$ otherwise.
\end{itemiz}
Finally, let us point out that most estimates hold for small enough
$\eps$ only, and often only for $\fP$-almost all $\omega\in\Omega$. We
will stress these facts only where confusion might arise.  

\begin{figure}
 \centerline{\psfig{figure=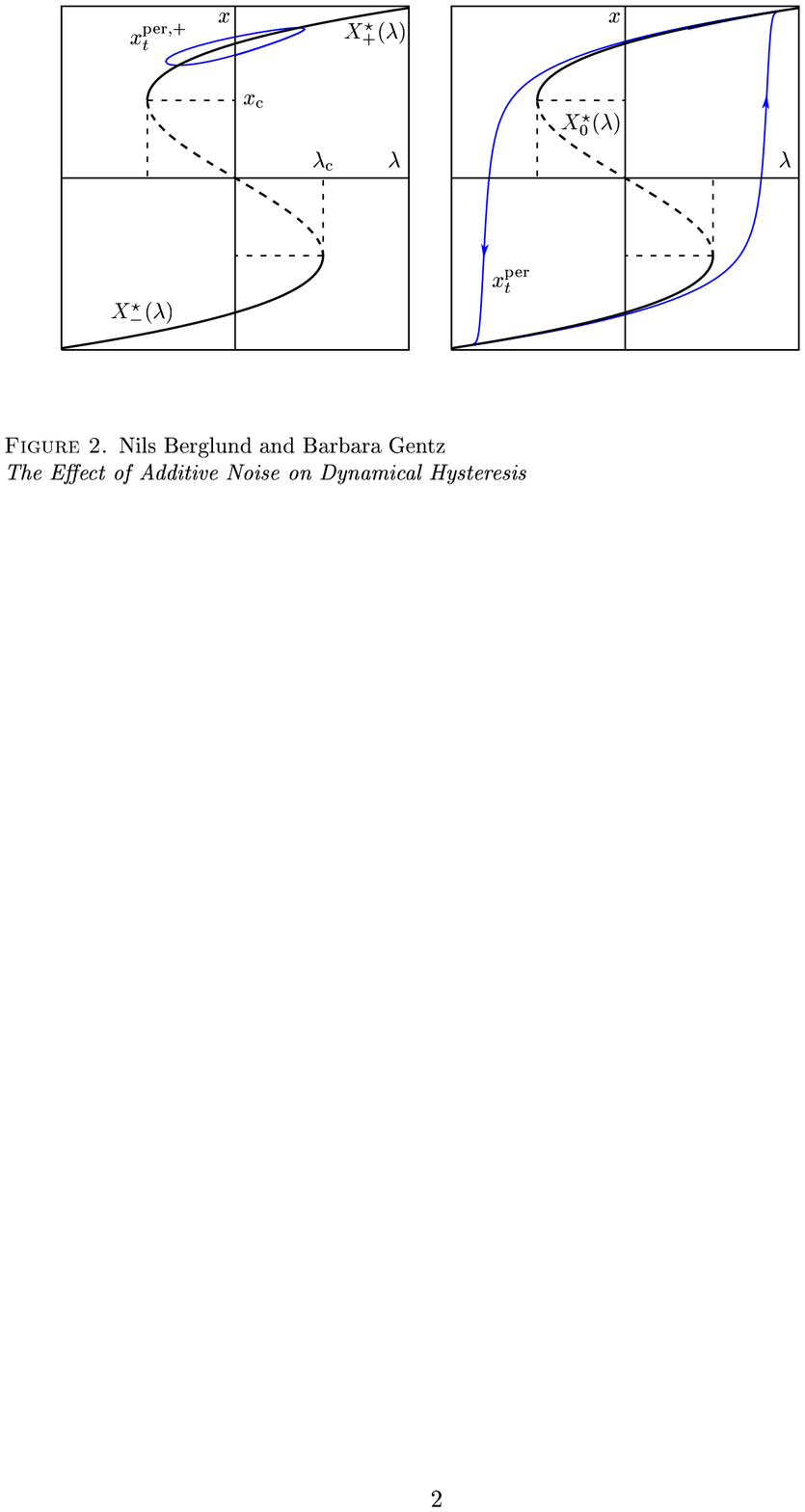,width=120mm,clip=t}}
 \caption[]
 {Equilibrium branches of $F$ (heavy curves) and periodic solutions of the
 deterministic equation (light curves), for $A<\lc$ (left) and $A>\lc$
 (right). For $A<\lc+\Order{\eps}$, the enclosed area is of order $\eps$
 while for $A>\lc+\Order{\eps}$, it is of order $\cA_0 +
 \eps^{2/3}(A-\lc)^{1/3}$.}
\label{fig2}
\end{figure}

\begin{figure}
 \centerline{\psfig{figure=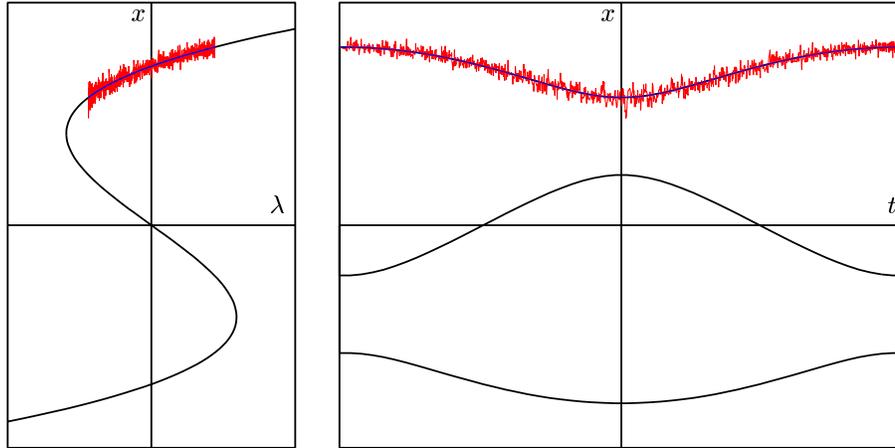,width=120mm,clip=t}}
 \caption[]
 {A sample path of Equation~\ref{SDE0} for $\eps=0.001$, $\sigma=0.05$ and
 $a_0=-0.1$, corresponding to the small amplitude regime.}
\label{fig3}
\end{figure}

\begin{figure}
 \centerline{\psfig{figure=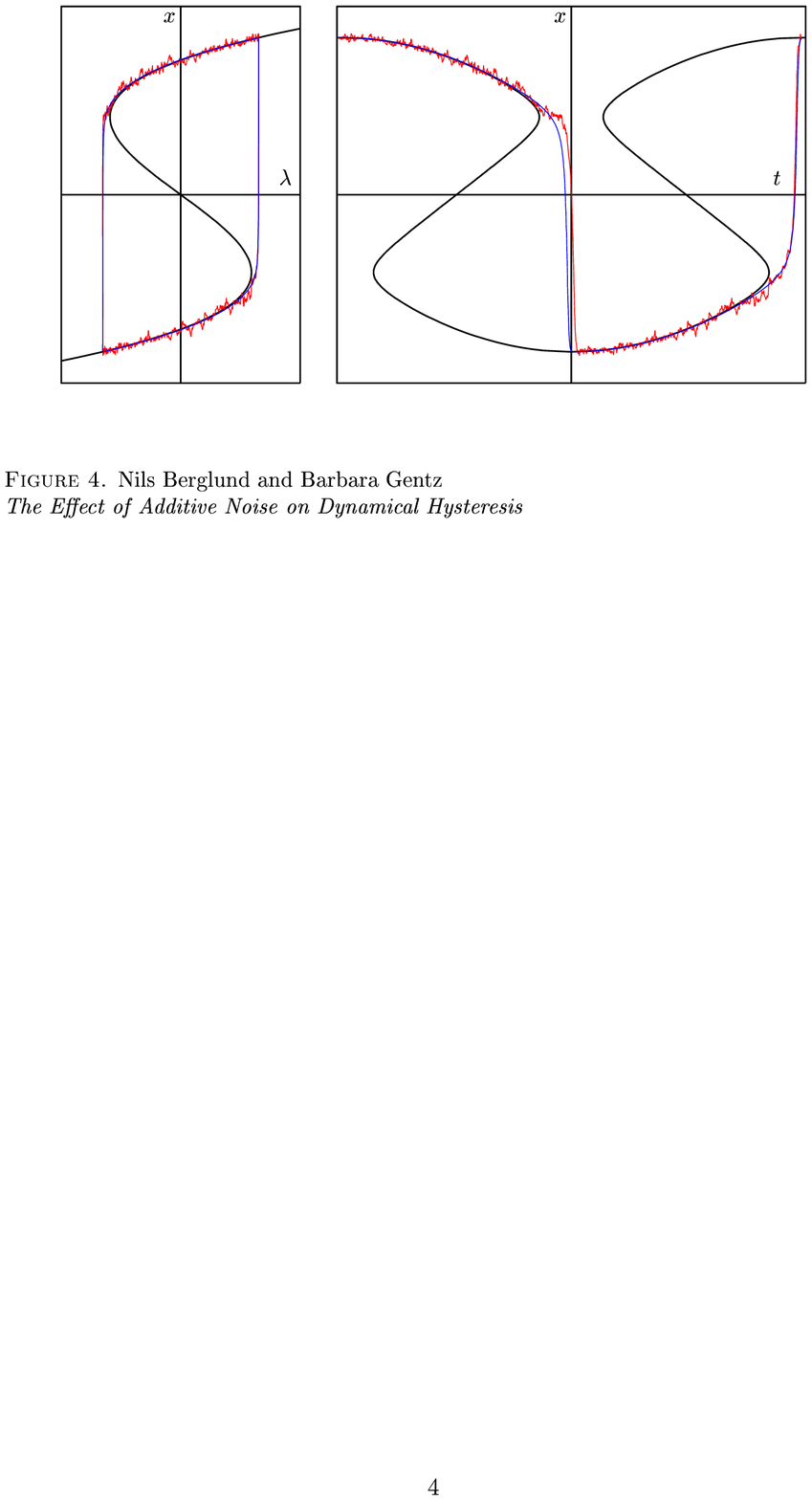,width=120mm,clip=t}}
 \caption[]
 {A sample path of Equation~\ref{SDE0} for $\eps=0.005$, $\sigma=0.04$ and
 $a_0=0.04$, corresponding to the large amplitude regime. A solution
 for $\sigma=0$ is shown for comparison.}
\label{fig4}
\end{figure}

\begin{figure}
 \centerline{\psfig{figure=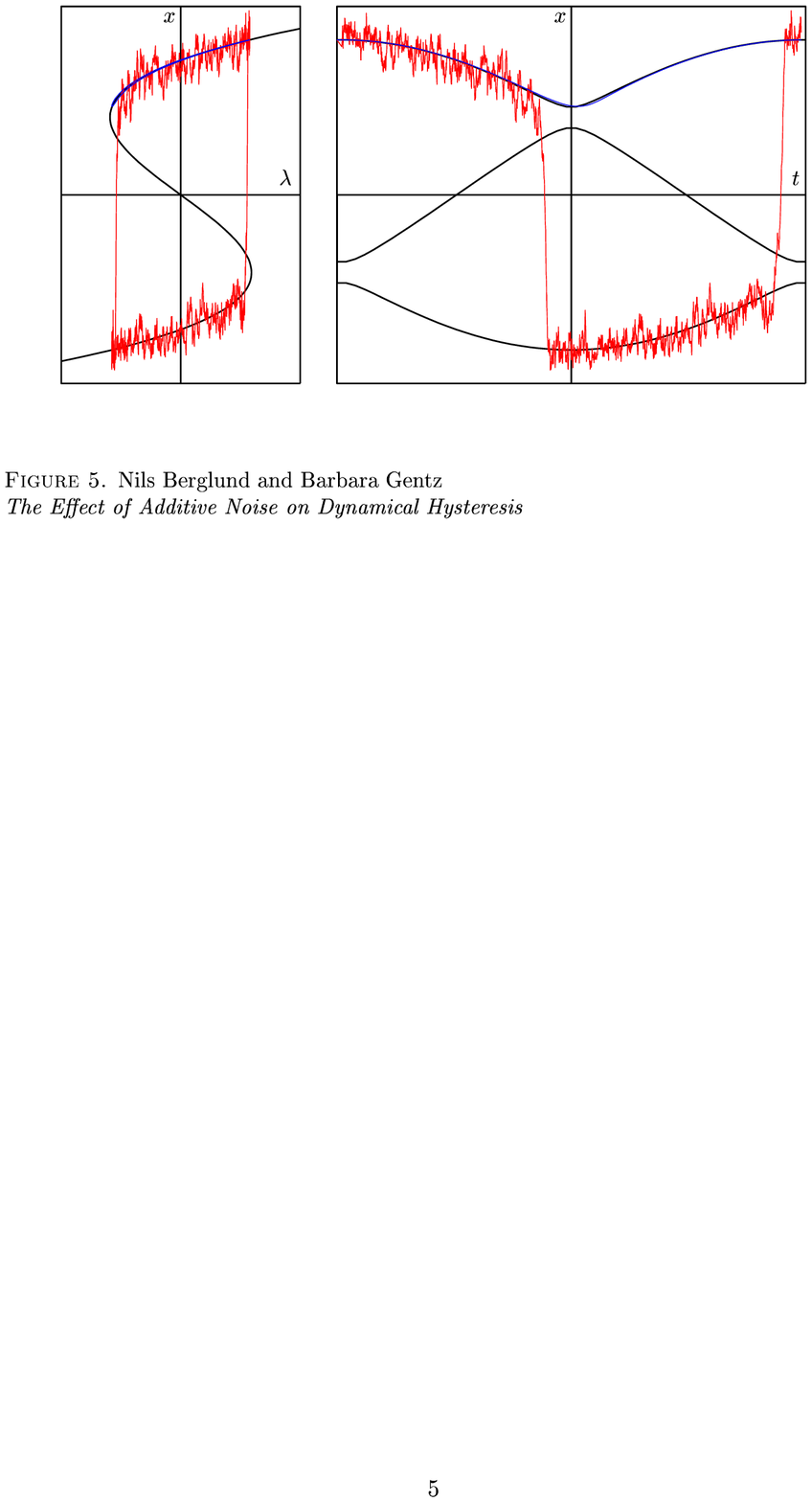,width=120mm,clip=t}}
 \caption[]
 {A sample path of Equation~\ref{SDE0} for $\eps=0.005$, $\sigma=0.16$ and
 $a_0=-0.01$, corresponding to the large noise regime. A solution for
 $\sigma=0$ is shown for comparison.}
\label{fig5}
\end{figure}

Let us first consider the deterministic case $\sigma=0$. It is convenient to
introduce the slow time $t = \eps s$, and rewrite \eqref{SDE0} for
$\sigma=0$ as 
\begin{equation}
\label{Det}
\eps\dtot {x_t}t = F(x_t,\lambda(t)).
\end{equation}
We start by discussing some properties of this equation, which will be
summarized in Theorem \ref{thm_det} below. As $\eps$ goes to zero,
solutions of \eqref{Det} are known to approach equilibrium branches of $F$,
that is, solutions of $F(x,\lambda)=0$ (see Figure~\ref{fig2}). Let $\lc =
2/(3\sqrt3)$. 
\begin{itemiz}
\item	
For $\abs{\lambda}<\lc$, $F$ has three equilibrium branches
$X^\star_-(\lambda) < X^\star_0(\lambda) < X^\star_+(\lambda)$, where
$X^\star_\pm(\lambda)$ are stable equilibria and $X^\star_0(\lambda)$ is an
unstable equilibrium of the associated frozen system $\dot x =
F(x,\lambda)$. 
\item	
At $\lambda=-\lc$, the branches $X^\star_+(\lambda)$ and
$X^\star_0(\lambda)$
undergo a saddle--node bifurcation, and $X^\star_+(-\lc) = X^\star_0(-\lc) =
\xc \defby 1/\sqrt3$. 
\item	
For $\lambda<-\lc$, $X^\star_-(\lambda)$ is the only equilibrium branch. 
\item	
A similar bifurcation occurs at $\lambda=+\lc$, where $X^\star_-(\lc)
= X^\star_0(\lc) = -\xc$. 
\item	
For $\lambda>\lc$, $X^\star_+(\lambda)$ is the only equilibrium branch. 
\end{itemiz}
We can thus expect a qualitative difference, in the limit $\eps\to0$,
between the regime $A<\lc$, where $F$ always derives from a double-well
potential, and the regime $A>\lc$, where $F$ has only one equilibrium part
of the time.

\begin{definition}
\label{def_hysteresis}
Let $\xperof{\eps}_t$ be a periodic solution of \eqref{Det}. We say that
this
solution \defwd{does not display hysteresis} if there exists a continuous
function $\lambda\mapsto X^\star(\lambda)$ such that 
\begin{equation}
\label{r2}
\lim_{\eps\to0} \xperof{\eps}_t = X^\star(\lambda(t)).
\end{equation}
If no such function exists, we say that $\xperof{\eps}_t$ \defwd{displays
hysteresis}. 
\end{definition}

If $A<\lc$, solutions starting near a stable equilibrium branch
$X^\star_+(\lambda)$ or $X^\star_-(\lambda)$ will remain close to that
branch, and relation \eqref{r2} holds with
$X^\star(\lambda)=X^\star_+(\lambda)$ or $X^\star_-(\lambda)$,
depending on the initial condition. If $A>\lc$, however, it turns out
that 
\begin{equation}
\label{r3}
\lim_{\eps\to0} \xperof{\eps}_t =
\begin{cases}
X^\star_+(\lambda(t)) 
&\text{if $\lambda(t)>\lc$ or if $\lambda(t)>-\lc$ and $\lambda'(t)<0$} \\
X^\star_-(\lambda(t))
&\text{otherwise.}
\end{cases}
\end{equation}
Thus the solution displays hysteresis since the instantaneous value of
$\lambda$ alone does not suffice to determine the state of the system in the
adiabatic limit. This so-called \defwd{hysteresis cycle\/} can be
characterized by its area, defined as 
\begin{equation}
\label{r4}
\cA(\eps) = -\int_{-1/2}^{1/2} \xperof{\eps}_t \lambda'(t)\6t.
\end{equation}
If $A<\lc$, we have $\lim_{\eps\to0}\cA(\eps)=0$, while for $A>\lc$, 
\begin{equation}
\label{r5}
\lim_{\eps\to0}\cA(\eps) = \cA_0 
\defby \int_{-\lc}^{\lc} (X^\star_+(\lambda) - X^\star_-(\lambda))\6\lambda
= \frac32.
\end{equation}
The situation is thus relatively simple in the limit $\eps\to0$. Since in
practice, however, the variation of $\lambda$ will not be infinitely
slow, it is important to understand what happens for small but
positive values of $\eps$. We summarize the necessary facts in the following
theorem. 

\begin{theorem}[Deterministic Case]
\label{thm_det}
There exist constants $\gamma_1>\gamma_0>0$ such that the following
behaviour holds for sufficiently small $\eps$.
\begin{itemiz}
\item	If $a_0 = A-\lc \leqs \gamma_0\eps$, Equation \eqref{Det} has
exactly two stable periodic solutions $\xperof{+}_t$ and $\xperof{-}_t$,
and one unstable periodic solution $\xperof{0}_t$. These solutions track,
respectively, the equilibrium branches $X^\star_\pm(\lambda(t))$ and
$X^\star_0(\lambda(t))$ at a distance not larger than
$\Order{\eps\abs{a_0}^{-1/2}\wedge\sqrt\eps\mskip1.5mu}$,  
and enclose an area
\begin{equation}
\label{r6}
\cA(\eps) \asymp \eps A.
\end{equation}
All solutions which do not start on $\xperof{0}_t$ are attracted
either by $\xperof{+}_t$ or by $\xperof{-}_t$.

\item	If $a_0 = A-\lc \geqs \gamma_1\eps$, Equation \eqref{Det}
admits exactly one periodic solution $\xper_t$. This solution is stable,
satisfies \eqref{r3} in the adiabatic limit, and encloses an area
$\cA(\eps)$ satisfying 
\begin{equation}
\label{r7}
\cA(\eps) - \cA_0 \asymp \eps^{2/3}a_0^{1/3}.
\end{equation}
\end{itemiz}
\end{theorem}

In the case where $a_0$ is of order $1$, the scaling law \eqref{r7} was
first obtained in \cite{Jung}. We outline the proof of Theorem~\ref{thm_det}
in Section~\ref{sec_det}. Note that in the transition zone
$\gamma_0\eps<a_0<\gamma_1\eps$, the situation is more complicated, since
more than two stable periodic orbits can coexist \cite{TO,BK}.

Let us now return to the stochastic differential equation \eqref{SDE0}. In
slow time $t=\eps s$, it can be written as
\begin{equation}
\label{SDE}
\6x_t = \frac1\eps F(x_t,\lambda(t)) \6t + \frac\sigma{\sqrt\eps}\6W_t.
\end{equation}
Let us fix, say, $t_0=-1/2$ as initial time, and some $x_0>0$ as initial
condition, such that the solution $\xdet_t$ of the deterministic equation
\eqref{Det} with the same initial condition is attracted by
$\xperof{+}_t$ or $\xper_t$, respectively. 
We denote by $x_t(\w)$ the solution of the SDE \eqref{SDE} with initial
condition $x_{t_0}=x_0$ for a given realization $\w$ of the Brownian motion,
and associate with it the area 
\begin{equation}
\label{r8}
\cA(\eps,\sigma;\w) = -\int_{-1/2}^{1/2} x_t(\w)\lambda'(t) \6t.
\end{equation}
Note that $\cA(\eps,\sigma;\w)$ also depends on $a_0=A-\lc$. We do not
stress this dependence here but consider $a_0$ as a (possibly constant)
function of $\eps$. Of course, since $x_t(\w)$ is not periodic in general,
the integral~\eqref{r8} does not represent the area of enclosed  by a
closed curve. However it is still physically meaningful since it describes
the energy dissipation if $x$ and $\lambda$ are thermodynamically conjugate
variables. One can check that for $\abs{x_0-\xperof{+}_{t_0}}$ sufficiently
small (but still of order one), $\abs{\xdet_t-\xperof{+}_t}$ decreases
exponentially fast in $(t-t_0)/\eps$ and thus $\cA(\eps,0)$ still behaves
like \eqref{r6}. The same is true for $\xper_t$ and the validity
of~\eqref{r7}.  

Our main purpose is to characterize the distribution of the random variable
$\cA(\eps,\sigma)$ as a function of the parameters $\eps$, $a_0=A-\lc$
and $\sigma$. The following three theorems describe the situation in three
different parameter regimes. 

\begin{theorem}[Case I -- Small amplitude regime]
\label{thm_smallamp}
Assume that $a_0=A-\lc\leqs\gamma_0\eps$ and that 
$\sigma\leqs(\abs{a_0}\vee\eps)^{3/4}$. Then there exist positive constants
$\kappa$, $h_1$, $h_2$, $c_0$ and $C$ such that the following properties
hold for sufficiently small $\eps$.
\begin{itemiz}
\item	The probability that a sample path starting near one potential
well crosses the potential barrier during one period is smaller than
\begin{equation}
\label{r9}
\frac C\eps \e^{-\kappa(\abs{a_0}\vee\eps)^{3/2}/\sigma^2}.
\end{equation}

\item	{\bf Case Ia:} Assume that either $a_0\geqs-\eps$ or
$\sigma\leqs\sqrt{\eps}/\abs{\log\abs{a_0}}$. Then the deviation of the
area $\cA(\eps,\sigma)$ from its deterministic value $\cA(\eps,0)$
satisfies 
\begin{equation}
\label{r10}
\bigprob{\abs{\cA(\eps,\sigma)-\cA(\eps,0)}\geqs H} \leqs 
\begin{cases}
\displaystyle\vrule height 12pt depth 14pt width 0pt
\frac C\eps \e^{-\kappa H^2/(\sigma^2\eps)}
&\text{for\/ $0\leqs H\leqs H_1(\eps,a_0)$,}
\\
\displaystyle\vrule height 14pt depth 10pt width 0pt
\frac C\eps \e^{-\kappa H^{4}/\sigma^2}
&\text{for\/ $H\geqs h_2$,}
\end{cases}
\end{equation}
where $H_1(\eps,a_0) = h_1\sqrt\eps(\abs{a_0}\vee\eps)^{3/4} \wedge
(\eps/\abs{\log(\abs{a_0}\vee\eps)})$. Furthermore, under the slightly
stronger assumption $\sigma\abs{\log\eps}\leqs
c_0(\abs{a_0}\vee\eps)^{3/4}$, 
\begin{align}
\label{r11a}
\bigabs{\bigexpec{\cA(\eps,\sigma)-\cA(\eps,0)}} &\leqs 
\Order{\sigma^2 \abs{\log(\abs{a_0}\vee\eps)}} \\  
\variance\bigset{\cA(\eps,\sigma)-\cA(\eps,0)}
&\asymp \sigma^2\eps.
\label{r11b}
\end{align}

\item	{\bf Case Ib:} Assume now that $a_0\leqs-\eps$ and
$\sigma\geqs\sqrt{\eps}/\abs{\log\abs{a_0}}$. Then \eqref{r10} still holds,
and in addition, we have 
\begin{equation}
\label{r12}
\bigprob{\abs{\cA(\eps,\sigma)-\cA(\eps,0)}\geqs H} \leqs 
\frac C\eps \e^{-\kappa H/(\sigma^2\abs{\log\abs{a_0}})}
\end{equation}
for $\eps/\abs{\log\abs{a_0}} \leqs H \leqs H_2(\eps,a_0)=
h_1\abs{a_0}^{3/2}\abs{\log\abs{a_0}}$. Moreover, under the slightly
stronger assumption $\sigma\abs{\log\eps}\leqs
c_0(\abs{a_0}\vee\eps)^{3/4}$, 
\begin{align}
\label{r13a}
\bigabs{\bigexpec{\cA(\eps,\sigma)-\cA(\eps,0)}} &\leqs 
\Order{\sigma^2 \abs{\log\abs{a_0}}} \\  
\variance\bigset{\cA(\eps,\sigma)-\cA(\eps,0)}
&\leqs \Order{\sigma^4 \abs{\log\abs{a_0}}^2}.
\label{r13b}
\end{align}
\end{itemiz}
\end{theorem}

The proof follows as a particular case of more general results presented in
Section~\ref{sec_near}. 

In Case Ia, the distribution of $\cA(\eps,\sigma)$ is close to a Gaussian
centred at $\cA(\eps,0)$. Both the expectation of
$\cA(\eps,\sigma)-\cA(\eps,0)$ and its standard deviation are smaller than
the deterministic value $\cA(\eps,0)\asymp\eps$. Thus one will still
observe, with a high probability, an area of the same order as the
deterministic one. 

In Case Ib, the distribution of $\cA(\eps,\sigma)$ becomes more spread
out, with a standard deviation possibly exceeding the deterministic
value $\cA(\eps,0)$. Thus, although typical values of
$\cA(\eps,\sigma)$ will still be small, the probability of negative
values is no longer negligible, and the deterministic scaling law
$\cA(\eps,0)\asymp\eps$ can no longer be observed. 

The quartic decay of the probability of deviations of order larger
than $1$ from the deterministic area is a consequence of the cubic
growth of the drift term $F$ for large $\abs{x}$. In fact, this
property holds in {\em all\/} parameter regimes, since it does not
depend on the details of the dynamics near the origin. For the sake of
brevity, we will not repeat this estimate in the other regimes. 

Note that there is a gap between $H\leqs H_1$ or $H_2$ and $H\geqs h_2$
where we do not describe the deviations. 
In fact, the distribution of $\cA(\eps,\sigma)$ will not be unimodal.
Sample paths are unlikely to jump from one potential well to the other
one, but if they do so, then  most likely near the instant of minimal
barrier height, producing a small peak in the distribution for areas
$\cA(\eps,\sigma)$ of order $1$.  

\begin{theorem}[Case II -- Large amplitude regime]
\label{thm_largeamp}
Assume that $a_0=A-\lc\geqs\gamma_1\eps$ and that 
$\sigma\leqs\esa0^{1/2}$. Then there exist positive constants $\kappa$,
$h_1$, $h_2$, $c_0$, $L_0$, $L_1$, $L_2$ and $C$ such that the
following properties hold for sufficiently small $\eps$.
\begin{itemiz}
\item	Let $\lambda^0$ denote the (random) value of $\lambda$ when $x_t$
changes sign for the first time. Then  
\begin{equation}
\label{r14}
\bigprob{\abs{\lambda^0}<\lc-L} \leqs \frac C\eps
\e^{-\kappa(L^{3/2}\vee\eps\sqrt{a_0})/\sigma^2}
\end{equation} 
for $-L_1\esa0^{2/3}\leqs L \leqs L_0/\abs{\log\esa0}$, and 
\begin{equation}
\label{r15}
\bigprob{\abs{\lambda^0}>\lc+L} \leqs 3 \,\indicator{(0,a_0]}(L)\,
\exp\biggset{-\frac\kappa{\sigma^2}\frac{L}{\esa0^{2/3}\abs{\log\esa0}}} 
\end{equation}
for $L\geqs L_2\esa0^{2/3}$.

\item	{\bf Case IIa:} Assume that $\sigma\leqs\esa0^{5/6}$. Then 
\begin{equation}
\label{r17}
\bigprob{\abs{\cA(\eps,\sigma)-\cA(\eps,0)}\geqs H} \leqs 
\frac C\eps \exp\biggset{-\kappa\frac{H^2}{\sigma^2\esa0^{1/3}}}
\qquad
\forall\,H\leqs h_1\eps\sqrt{a_0}.
\end{equation} 
Furthermore, under the slightly stronger assumption
$\sigma\abs{\log\eps}\leqs c_0\esa0^{5/6}$, 
\begin{align}
\label{r18a}
\bigabs{\bigexpec{\cA(\eps,\sigma)-\cA(\eps,0)}} & \leqs 
\biggOrder{\frac{\sigma^2\abs{\log\eps}}{\esa0^{2/3}}} \\  
\variance\bigset{\cA(\eps,\sigma)-\cA(\eps,0)}
&\asymp \sigma^2\esa0^{1/3}.
\label{r18b}
\end{align}

\item	{\bf Case IIb:} Assume now that
$\esa0^{5/6}\leqs\sigma\leqs\esa0^{1/2}$. Then 
\begin{align}
\label{r19a}
\bigprob{\cA(\eps,\sigma)-\cA(\eps,0)\leqs -H} &\leqs 
\frac C\eps \e^{-\kappa H^{3/2}/\sigma^2} \\
\bigprob{\cA(\eps,\sigma)-\cA(\eps,0)\geqs +H} &\leqs 
\frac C\eps \exp\biggset{-\kappa\frac{\esa0^{1/3}H}{\sigma^2\abs{\log\esa0}}}
\label{r19b}
\end{align}
for $h_1\esa0^{2/3}\abs{\log\esa0}\leqs H\leqs
h_2\esa0^{1/3}\abs{\log(\eps^{2/3}a_0^{-1/6})}$.  As a consequence,  if
$\sigma \leqs c_0\esa0 / \sqrt{\abs{\log\eps}}$, then  expectation and
standard deviation of $\cA(\eps,\sigma)-\cA(\eps,0)$ are both at most of
order $\esa0^{2/3}\abs{\log\esa0}$. 
\end{itemiz}
\end{theorem}

The proof is given in Section~\ref{sec_lac}. 

The estimates \eqref{r14} and \eqref{r15}  show that the value of the
parameter (e.\,g.\ the magnetic field) for which $x_t$ changes sign is most
likely $\pm(\lc+\Order{\esa0^{2/3}})$, which corresponds to the
deterministic value. 

In Case IIa, the distribution of $\cA(\eps,\sigma)$ is again close to a
Gaussian in a neighbourhood of $\cA(\eps,0)$. Both the expectation of
$\cA(\eps,\sigma)-\cA(\eps,0)$ and its standard deviation are smaller than
the deterministic value of $\cA(\eps,0)-\cA_0\asymp\esa0^{2/3}$. Thus one
will still observe, with a high probability, an area of the same order as
the deterministic one. 

In Case IIb, we can only show that $\cA(\eps,\sigma)$ is likely to belong
to an interval of size $\esa0^{2/3}\abs{\log\esa0}$ centred at the
deterministic value, so that $\cA(\eps,\sigma)-\cA_0$ is not necessarily
positive with probability close to $1$. There is a gap between the
estimates \eqref{r19a} and \eqref{r19b} outside this interval, and the
trivial bound $1$ inside the interval. This is due to the large spreading
of paths during the jump. However, this result may conceivably fall short
of being optimal. 

\begin{theorem}[Case III -- Large noise regime]
\label{thm_largenoise}
Assume that either $a_0\leqs\eps$ and $\sigma>(\abs{a_0}\vee\eps)^{3/4}$, or
$a_0\geqs\eps$ and $\sigma > \esa0^{1/2}$. Then there exists a
(deterministic) reference area $\hat\cA$, satisfying 
\begin{equation}
\label{r21}
\hat\cA - \cA_0 \asymp -\sigma^{4/3},
\end{equation}
and positive constants $\kappa$, $h_1$, $h_2$, $c_0$, $c_1$, $c_2$ and
$C$, such that the following properties hold for sufficiently small $\eps$.
\begin{itemiz}
\item	{\bf Case IIIa:} Either $a_0\leqs\eps$ or $\sigma>a_0^{3/4}$. Then
the deviation of the area $\cA(\eps,\sigma)$ from the reference value
$\hat\cA$ satisfies 
\begin{align}
\label{r22a}
\bigprob{\cA(\eps,\sigma)-\hat\cA < -H} &\leqs 
\frac C\eps \e^{-\kappa H^{3/2}/\sigma^2} + 
\frac32 \e^{-\kappa\sigma^{4/3}/(\eps\abs{\log\sigma})} \\
\bigprob{\cA(\eps,\sigma)-\hat\cA > +H} &\leqs 
\frac C\eps \e^{-\kappa H/(\sigma^2\abs{\log\sigma})} + 
\frac32 \indicator{[0,h_2\sigma^{4/3})}(H) \e^{-\kappa
H/(\eps\abs{\log\sigma})}
\label{r22b}
\end{align}
for $0\leqs H\leqs h_1\sigma^{2/3}\abs{\log\sigma}$. Moreover, if
the noise intensity satisfies
$c_1\eps<\sigma^{4/3}/\abs{\log\sigma}^2$ and
$\sigma^{2/3}\abs{\log\sigma} \leqs c_2/\abs{\log\eps}$, then 
\begin{align}
\label{r23a}
\bigexpec{\cA(\eps,\sigma)-\hat\cA} 
&\in [-C \sigma^{4/3} \abs{\log\eps}^{2/3},
C (\eps \vee \sigma^2\abs{\log\eps})\abs{\log\sigma}] \\
\label{r23b}
\variance\bigset{\cA(\eps,\sigma)-\hat\cA} 
&\leqs C \bigpar{\sigma^{4/3} \abs{\log\eps}^{2/3}}^2. 
\end{align}

\item	{\bf Case IIIb:} $a_0\geqs\eps$ and $\sigma\leqs a_0^{3/4}$. 
Let $\ell_0=\abs{\log(\sigma^{4/3}/\sqrt{a_0})}$. Then
\begin{align}
\label{r24a}
\bigprob{\cA(\eps,\sigma)-\hat\cA < -H} &\leqs 
\frac C\eps \e^{-\kappa H^{3/2}/\sigma^2} + 
\frac32 \e^{-\kappa\sigma^2/(\eps\sqrt{a_0}\abs{\log\sigma})} \\
\bigprob{\cA(\eps,\sigma)-\hat\cA > +H} &\leqs 
\frac C\eps \e^{-\kappa H/(\sigma^2\ell_0)} + 
\frac32 \indicator{[0,h_2a_0)}(H) \e^{-\kappa\sigma^{2/3}
H/(\eps\sqrt{a_0}\abs{\log\sigma})}
\label{r24b}
\end{align}
for $0\leqs H\leqs h_1\sigma^{2/3}\ell_0$. In addition, if  $\sigma \geqs
c_1\esa0^{1/2} \abs{\log\eps}$, then
\begin{align}
\label{r25a}
\bigexpec{\cA(\eps,\sigma)-\hat\cA} 
&\in [-C \sigma^{4/3} \abs{\log\eps}^{2/3},
C (\sigma^2\ell_0\abs{\log\eps} \vee
\eps\sqrt{a_0}\abs{\log\sigma}/\sigma^{2/3})] \\
\label{r25b}
\variance\bigset{\cA(\eps,\sigma)-\hat\cA} 
&\leqs C \bigpar{\sigma^{4/3} \abs{\log\eps}^{2/3}}^2. 
\end{align}
\end{itemiz}
The value $\lambda^0$ of $\lambda$ at the first time $x_t$ reaches $0$
behaves in a similar way as the area, when compared to a reference value
$\hat\lambda$ equal to $\lc-\Order{\sigma^{4/3}}$. 
\end{theorem}

The proof is given in Section~\ref{sec_lnr}. 

The main feature in this parameter regime is that the noise intensity is
sufficiently large to drive $x_t$ over the potential barrier before it
reaches its minimal height or even vanishes. The barrier is typically crossed
when $\abs{\lambda}$ equals $\lc-\Order{\sigma^{4/3}}$. 

The distribution of $\cA(\eps,\sigma)$ decays faster to the right of
$\hat\cA$ than to the left. The probability that $\cA(\eps,\sigma)$
exceeds $\cA_0$ is very small (unless $\sigma$ approaches its
threshold value), so that it is indeed likely to observe an area that
is smaller than $\cA_0$, by an amount of order $\sigma^{4/3}$. 


\section{Deterministic case}
\label{sec_det}

In this section we discuss the deterministic equation
\begin{equation}
\label{det1}
\eps\dtot xt = x - x^3 + \lambda(t)
\end{equation}
with initial condition $x_{-1/2}=x_0$. Recall that we are interested
in the case $\lambda(t) = -A \cos(2\pi t)$ with $A=\lc+a_0$. Since
this equation has already been studied in \cite{Jung,TO,BK,BG2}, we
only outline the main properties without proofs. 


\subsection{The case $a_0\leqs\gamma_0\eps$}
\label{ssec_det1}

The simplest situation occurs when $a_0$ is negative and of order $1$. Then
the three curves $x^\star_\pm(t)=X^\star_\pm(\lambda(t))$ and
$x^\star_0(t)=X^\star_0(\lambda(t))$ are uniformly hyperbolic equilibrium
curves of the associated family of frozen systems $\dot x=x-x^3+\lambda$.
Thus Tihonov's theorem \cite{Tihonov,Grad} shows the existence of
particular solutions $x^\pm_t$ and $x^0_t$ tracking, respectively,
$x^\star_\pm(t)$ and $x^\star_0(t)$ at a distance of order $\eps$. These
solutions are not necessarily periodic, but the curves $x^\pm_t$ attract a
\nbh\ of order $1$ exponentially fast. Thus the Poincar\'e map
$P:x_{-1/2}\mapsto x_{1/2}$ maps \nbh s of order $1$ of
$x^\star_\pm(-1/2)$, respectively, to two exponentially small intervals
containing $\smash{x^\pm_{1/2}}$. This implies the existence of a unique
fixed point (corresponding to a periodic orbit) in each interval. A similar
statement is true for $P^{-1}$ in a \nbh\ of $x^\star_0(1/2)$. The fact
that $x_t$ is monotonous between the equilibrium branches excludes the
existence of other periodic orbits (see,
e.g.~\cite[Proposition~4.8]{Thesis}). 

If $a_0$ is a small parameter, Tihonov's theorem can also be applied
outside a given interval $[-T,T]$, $T$ a constant of order $1$, while the
dynamics in $[-T,T]$ has to be analysed separately. For $\abs{t}\geqs
c_0(\abs{a_0}\vee\eps)^{1/2}$, $c_0$ a sufficiently large
constant, one can consider the deviation $y_t=x_t-x^\star_+(t)$, which obeys
the equation
\begin{equation}
\label{det2}
\eps\dtot yt = a^\star_+(t)y + b^\star_+(y,t) - \eps\dtot{x^\star_+}t,
\end{equation}
where
\begin{equation}
\label{det3}
\begin{split}
a^\star_+(t) &= 1-3(x^\star_+(t))^2 \\
b^\star_+(y,t) &= -y^2\bigbrak{3x^\star_+(t)+y}.
\end{split}
\end{equation}
This equation is used to show (see \cite[Section~4.1]{BG2}) that 
\begin{equation}
\label{det4}
y_t \asymp \frac\eps{\abs t} 
\qquad\qquad
\text{for $-T\leqs t \leqs -c_0(\abs{a_0}\vee\eps)^{1/2}$.}
\end{equation}
In the case $a_0\leqs -\eps/\gamma_0$ for some small enough
$\gamma_0>0$, one also obtains from \eqref{det2} that for
$\abs{t}\leqs c_0\sqrt{\abs{a_0}}$,  
\begin{equation}
\label{det5}
y_t = -t\mskip1.5mu C_1(t) + C_2(t) 
\qquad\text{with}\qquad 
C_1(t)\asymp\frac{\eps}{\abs{a_0}},
\quad 
C_2(t)\asymp\frac{\eps^2}{\abs{a_0}^{3/2}},
\end{equation}
which implies that $y_t$ becomes negative at a time of order
$\eps\abs{a_0}^{-1/2}$. If $-\eps/\gamma_0\leqs a_0\leqs\gamma_0\eps$,
the dynamics for $\abs{t}\leqs c_0\sqrt\eps$ is analysed by the change
of variables 
\begin{equation} 
\label{det6}
t = c\sqrt\eps s, 
\qquad\qquad
x = \frac1{\sqrt3} \Bigpar{1+\frac{\sqrt\eps}c z},
\end{equation}
where $2\pi^2Ac^4=1$. Then $z$ obeys a perturbation of order $\sqrt\eps$ of
the Riccati equation
\begin{equation}
\label{det7}
\dtot zs = s^2 - z^2 - \delta,
\qquad\qquad
\text{where $\delta = \displaystyle \sqrt3 \,c^2 \,\frac{a_0}\eps \leqs
\sqrt3 \,c^2 \, \gamma_0$,}
\end{equation}
which can be used to show that for $\gamma_0$ small enough, $z_s\asymp 1$
for $s$ of order $1$. It follows that for {\em all\/} $a_0\leqs
\gamma_0\eps$,  
\begin{equation}
\label{det8}
x_t - \frac1{\sqrt3} \asymp (\abs{a_0}\vee\eps)^{1/2} 
\qquad\qquad
\text{for $\abs{t}\leqs c_0(\abs{a_0}\vee\eps)^{1/2}$.}
\end{equation}
Finally, one obtains as before that 
\begin{equation}
\label{det9}
y_t \asymp -\frac\eps{\abs t} 
\qquad\qquad
\text{for $c_0(\abs{a_0}\vee\eps)^{1/2}\leqs t\leqs T$.}
\end{equation}
Hence there is a solution of \eqref{det1} tracking $x^\star_+(t)$ at a
distance of order
$\eps/(\abs{t}\vee\sqrt{\abs{a_0}}\vee\sqrt\eps\mskip1.5mu)$ (if $a_0>0$,
$x^\star_+(t)$ does not exist during a small time interval, but this gap is
too small for $x_t$ to slip through). Later we will use the fact that the
linearization of $F$ around this solution satisfies
\begin{equation}
\label{det10}
a(t) \defby \dpar{}x \bigbrak{x - x^3 - A \cos(2\pi t)} \Bigevalat{x=x_t} 
\asymp -(\abs{t}\vee\sqrt{\abs{a_0}}\vee\sqrt\eps\mskip1.5mu).
\end{equation}
Furthermore, we will need that fact that $\abs{a'(t)}$ is bounded above by a
constant independent of $\eps$ and $a_0$. This can be shown by using the
relation
\begin{equation}
\label{det10a}
a'(t) = -6x_t\dtot{}t x_t 
= -6 x_t \frac1\eps F(x_t,\lambda(t)). 
\end{equation}
The cases $\abs{t}\geqs (\abs{a_0}\vee\eps)^{1/2}$ and $a_0\leqs
-\eps/\gamma_0$ can be treated by expanding $F$ around the equilibrium branch
$x^\star_+(t)$ that $x_t$ is tracking, and using the
estimates~\eqref{det4}, \eqref{det5}, \eqref{det9} for
$x_t-x^\star_+(t)$. The remaining case can be treated by 
considering \eqref{det7} directly.  

In addition, these estimates show that
\begin{equation}
\label{det11}
\cA(\eps,0) = 2\pi A \int_{-1/2}^{1/2} x_t \sin(2\pi t)\6t \asymp \eps A.
\end{equation}

Similar properties hold for solutions tracking $x^\star_-(t)$ and
$x^\star_0(t)$, and the above arguments on the Poincar\'e map can be
repeated to show the existence of two stable and one unstable periodic
orbit.   


\subsection{The case $a_0\geqs\gamma_1\eps$}
\label{ssec_det2}

Let $\tc$ be the solution of $A\cos(2\pi\tc)=\lc$ in $[0,1/4]$. The
equilibrium branches $x^\star_+$ and $x^\star_0$ bifurcate at the point
$(-\tc,\xc)$, where $\xc=1/\sqrt3$. The translation 
\begin{equation}
\label{det12}
t = -\tc+s, 
\qquad\qquad
x = \xc+y
\end{equation}
yields the equation
\begin{equation}
\label{det13}
\eps\dtot ys = \mu(s) - \sqrt3 y^2 - y^3,
\end{equation}
where
\begin{align}
\nonumber
\mu(s) 
&= \lc - A\cos(2\pi(-\tc+s)) \\
\nonumber
&= \lc(1-\cos(2\pi s)) - \sqrt{A^2-\lc^2}\, \sin(2\pi s) \\
&= -2\pi \sqrt{2a_0\lc+a_0^2}\, s + 2\pi^2 \lc s^2 + \Order{s^3}.
\label{det14}
\end{align}
As before, one shows that $y_s$ tracks the equilibrium branch
$y^\star_+(s) = x^\star_+(-\tc+s)-\xc$ at a distance scaling like
$\eps/\abs{s}$ for $s\leqs-\eps^{2/3}\smash{a_0^{-1/6}}$. For larger
times, we use the scaling 
\begin{equation}
\label{det15}
s = c\eps^{2/3}a_0^{-1/6} u,
\qquad\qquad
y = \frac1{\sqrt3 c} \eps^{1/3} a_0^{1/6} z
\end{equation}
which yields, for an appropriate choice of $c$, a perturbation of order
$\eps^{1/3}a_0^{1/6}$ of the Riccati equation
\begin{equation}
\label{det16}
\dtot zu = \hat\mu(u) - z^2, 
\qquad\qquad
\text{where $\hat\mu(u) = -u + \Order{\eps^{2/3}a_0^{-2/3} u^2}$.}
\end{equation}
For sufficiently large $\gamma_1$ (recall that $a_0\geqs\gamma_1\eps$), one
can show that $z_u$ reaches a value of order $-1$ in a time of order $1$,
and thus $y_s$ reaches order $-\eps^{1/3}\smash{a_0^{1/6}}$ for some $s$ of
order $\eps^{2/3}\smash{a_0^{-1/6}}$. Finally, the fact that the right-hand
side of \eqref{det13} is smaller than $-y^2$ for sufficiently small $s$ and
$y$ can be used to show that $y_s$ reaches values of order $-1$ after
another time of order $\eps^{2/3}\smash{a_0^{-1/6}}$. For larger
times, $y_s$ is quickly attracted by the lower stable equilibrium
branch. Later we will use the fact that the linearization around $y_s$
satisfies 
\begin{equation}
\label{det17}
a(-\tc+s) \defby \dpar{}y \bigbrak{\mu(s) - \sqrt3 y^2 - y^3}
\Bigevalat{y=y_s}  
\asymp -(\abs{s}\vee a_0^{1/4}\sqrt{\abs{s}}\vee\eps^{1/3}a_0^{1/6})
\end{equation}
for $s \leqs c_0\eps^{2/3}a_0^{-1/6}$.
Furthermore, one can check that $\abs{a'(t)}$ is bounded above by a
constant times $(a_0/\eps)^{1/3}$ for $t\leqs -\tc +
\Order{\smash{\eps^{2/3}a_0^{-1/6}}}$.  It is easy to see that the solution
we constructed encloses an area satisfying $\cA-\cA_0\asymp
\eps^{2/3}\smash{a_0^{1/3}}$, where the main contribution comes from the
delayed jump from a \nbh\ of $x^\star_+(t)$ to a \nbh\ of $x^\star_-(t)$.
Since $x^\star_+(t)$ is the only equilibrium branch near $t_0=-1/2$, the
Poincar\'e map contracts any interval of order $1$ containing
$x^\star_+(t)$ to an exponentially small \nbh\ of $x_{1/2}$, which implies
the existence of a unique periodic orbit.  


\section{The random motion near stable equilibrium branches}
\label{sec_near}

We consider now the SDE 
\begin{equation}
\label{nse1}
\6x_t = \frac1\eps \bigbrak{x_t-x_t^3-A\cos(2\pi t)} \6t +
\frac\sigma{\sqrt\eps}\6W_t
\end{equation}
with a given (deterministic) initial condition $x_{t_0}=x_0$. Let $\xdet_t$
denote the solution of the deterministic equation \eqref{det1} with the
same initial condition. We will start by investigating the difference
$x_t-\xdet_t$, and then derive some properties of the area delimited by
this difference.


\subsection{Noise-induced deviations from the deterministic solution}
\label{ssec_soldev}

The difference $y_t=x_t-\xdet_t$ obeys the SDE
\begin{equation}
\label{nse2}
\6y_t = \frac1\eps \bigbrak{a(t) y_t + b(y_t,t)}\6t +
\frac\sigma{\sqrt\eps}\6W_t,
\qquad
y_{t_0}=0,
\end{equation}
where
\begin{equation}
\label{nse3}
\begin{split}
a(t) &= 1-3(\xdet_t)^2 \\
b(y,t) &= -y^2 \bigbrak{3\xdet_t+y}.
\end{split}
\end{equation}
In this section, we are interested in situations where $\xdet_s$ is
attracting up to time $t$, that is, we assume that $a(s)<0$ for
$t_0\leqs s\leqs t$. Results from the previous section (see \eqref{det10} and
\eqref{det17}) show that this is true when the following condition is
satisfied. 
\begin{assump}[Stable case]
\label{assump_stable}
Assume 
\begin{itemiz}
\item	either $a_0\leqs\gamma_0\eps$ and $t$ arbitrary
\item   or $a_0>\gamma_1\eps$, $t_0\geqs-1/2$ and $t\leqs
t^\star\defby -\tc+c_0\eps^{2/3}a_0^{-1/6}$. 
\end{itemiz}
\end{assump}

\noindent
For the sake of brevity, we will refer to these assumptions as {\it
stable case}. 

If we were to omit the nonlinear term $b(y,t)$ in \eqref{nse2}, the
solution $y_t$ of the equation would be normally distributed with mean
zero and variance
\begin{equation}
\label{nse4}
v(t) = \frac{\sigma^2}\eps \int_{t_0}^t \e^{2\alpha(t,s)/\eps}\6s, 
\qquad\qquad
\text{where $\alpha(t,s)=\int_s^t a(u)\6u$.}
\end{equation}
It is straightforward to show that the function 
\begin{equation}
\label{nse5}
\z(t) \defby \frac1{2\abs{a(t_0)}}\e^{2\alpha(t,t_0)/\eps} 
+ \frac1\eps \int_{t_0}^t \e^{2\alpha(t,s)/\eps}\6s 
\end{equation}
satisfies, in both cases summarized in Assumption~\ref{assump_stable},
\begin{equation}
\label{nse6}
\z(t)\asymp \frac1{2\abs{a(t)}}. 
\end{equation}
Thus in the linear case, the standard deviation of $y_t$ is smaller than
$\sigma\sqrt{\z(t)}$. The following result applies to the whole path
$\set{y_s}_{t_0\leqs s\leqs t}$ of the nonlinear equation \eqref{nse2}.

\begin{prop}
\label{prop_nse1}
Let 
\begin{equation}
\label{nse7}
\hat\z(t) = \sup_{t_0\leqs s\leqs t} \z(s).
\end{equation}
In the stable case, there exists a constant $h_0$ ($h_0\asymp 1$) such
that for all $h\leqs h_0\hat\z(t)^{-3/2}$, 
\begin{equation}
\label{nse8}
\biggprobin{t_0,0}{\sup_{t_0\leqs s\leqs t} 
\frac{\abs{y_s}}{\sqrt{\z(s)}} > h} 
\leqs 
\biggpar{\frac{\abs{\alpha(t,t_0)}}{\eps^2} + 2} 
\exp\biggset{-\frac12 \frac{h^2}{\sigma^2}
\Bigpar{1-\Order{\eps}-\bigOrder{h\hat\z(t)^{3/2}}}}.
\end{equation}
\end{prop}
\nobreak
\begin{proof}
The case $a_0\leqs0$ is the one considered in~\cite[Theorem~2.6]{BG2},
and the other cases can be proved in exactly the same way. 
\end{proof}
\goodbreak

If we do not care for the precise value of the exponent in~\eqref{nse8}, an
obvious modification in the proof yields the bound
\begin{equation}
\label{nse9}
\biggprobin{t_0,0}{\sup_{t_0\leqs s\leqs t} 
\frac{\abs{y_s}}{\sqrt{\z(s)}} > h} 
\leqs 
C\biggpar{\frac{t-t_0}\eps + 1} \e^{-\kappa h^2/\sigma^2}, 
\end{equation}
for all $h\leqs h_0\hat\z(t)^{-3/2}$, where $C$ and $\kappa$ are
positive constants. This estimate shows that in 
the time interval $[t_0,t]$, the typical spreading of paths is of order
$\sigma\sqrt{\z(s)}$ for $\sigma\ll\hat\z(t)^{-3/2}$. It allows to bound
the probability of deviations up to order $\hat\z(t)^{-1}$. On the other
hand, the special (cubic) form of the drift term in Equation \eqref{nse2}
allows for a bound on deviations of order larger than~$1$.

\begin{prop}
\label{prop_nse2}
There exist constants $L_0, C, \kappa>0$ such that for all $L\geqs L_0$ and
all $y_0\leqs L_0/2$, 
\begin{equation}
\label{nse10}
\biggprobin{t_0,y_0}{\sup_{t_0\leqs s\leqs t} y_s > L} 
\leqs C\biggpar{\frac{t-t_0}\eps + 1} \e^{-\kappa L^4/\sigma^2}.
\end{equation}
\end{prop}
\nobreak
\begin{proof}
First note that we are working in slow time. The estimate is
classical for $t-t_0\leqs \eps$, and starting from there, \eqref{nse10}
can be obtained by considering a partition of the interval $[t_0,t]$
with spacing proportional to $\eps/(t-t_0)$.
\end{proof}

\begin{remark}
\label{cubic_drift}
Note that the preceding proposition holds for all $a_0$, but $L_0$ may
depend on the amplitude $A$.
\end{remark}

The following proposition gives bounds on the moments of $y_t$. These
bounds hold whenever the estimates~\eqref{nse9} and~\eqref{nse10} are
satisfied, and we do not need to assume that $a(s)<0$ holds for all $s$.

\begin{prop}
\label{cor_nse_gen}
Fix $t>t_0$ such that $t-t_0$ is at most of order $1$ and assume that there
exists an $h_0(\eps,t)>0$ such that~\eqref{nse9} holds for 
all $h\leqs h_0(\eps,t)$. Then there exist constants $K,M>0$ such that 
\begin{align}
\label{nse12a_gen}
\bigexpecin{t_0,0}{\abs{y_t}^{2k}} &\leqs k!\,M^k\sigma^{2k}\z(t)^k
\biggpar{1+\frac{\abs{\log\eps}^k}{k!}}, 
\qquad k\in\N,
\end{align}
whenever $\sigma$ satisfies $K\sigma\sqrt{\abs{\log(\eps\sigma^2)}}
\leqs h_0(\eps,t)$. 
\end{prop}
\begin{proof}
We will only prove the case $k=1$, as the general case follows along
the same lines. 
Let $\gamma=K\sigma\sqrt{\abs{\log\eps}}$ for some constant
$K>0$ to be chosen later and set $h=h_0(\eps,t)$. Note
that, under our condition on $\sigma$, we may assume $\gamma<h$. 
Let $\z_0=\inf_{s\in[t_0,t]}\z(s)$. We write the expectation of
$y_t^2$ as $E_1+E$, where 
\begin{equation}
\label{cor_nse:2new}
E_1 = \biggexpecin{t_0,0}{y_t^2 \,\bigindexfct{\sup_{t_0\leqs s\leqs
t}\frac{\abs{y_s}}{\sqrt{\z(s)}}\leqs\gamma}}
\quad\text{\ and\ }\quad
E = \biggexpecin{t_0,0}{y_t^2 \,\bigindexfct{\sup_{t_0\leqs s\leqs
t}\frac{\abs{y_s}}{\sqrt{\z(s)}}>\gamma}}. 
\end{equation}
The first term $E_1$ can be estimated trivially, namely by $E_1 \leqs
\gamma^2\z(t)$. To estimate the second term $E$, we employ
integration by parts, thereby obtaining
\begin{equation}
\label{cor_nse:4new} 
E \leqs \z(t) \int_0^\infty 2z \, \biggprobin{t_0,0}{\sup_{t_0\leqs s\leqs
t}\frac{\abs{y_s}}{\sqrt{\z(s)}} \geqs\gamma \vee z}\6z.
\end{equation}
We now split the integral at $\gamma$, $h$ and $L_0/\sqrt{\z_0}$ and
estimate the resulting terms separately. By~\eqref{nse9},
\begin{align}
\nonumber
E_2 & =  \z(t) \int_0^{L_0/\sqrt{\z_0}} 2z \,
\biggprobin{t_0,0}{\sup_{t_0\leqs s\leqs t}
\frac{\abs{y_s}}{\sqrt{\z(s)}} \geqs\gamma \vee z}\6z \\ 
\label{cor_nse:5}
    & \leqs \z(t)  C\Bigpar{\frac{t-t_0}\eps + 1}
\biggbrak{\Bigpar{\gamma^2
+\frac{\sigma^2}{\kappa}}\e^{-\kappa\gamma^2/\sigma^2} +
\frac{L_0^2}{\z_0} \e^{-\kappa h^2/\sigma^2}}. 
\end{align}
Estimating the remaining part of the integral with the help
of~\eqref{nse10}, 
\begin{equation}
\label{cor_nse:5b}
E_3 =  \z(t) \int_{L_0/\sqrt{\z_0}}^\infty 2z \,
\biggprobin{t_0,0}{\sup_{t_0\leqs s\leqs t}
\frac{\abs{y_s}}{\sqrt{\z(s)}} \geqs\gamma \vee z}\6z 
\leqs \z(t)  C\Bigpar{\frac{t-t_0}\eps + 1}
\frac{\sigma^2}{2\kappa L_0^2\z_0}\e^{-\kappa L_0^4/\sigma^2}
\end{equation}
follows. Since the expectation of $y_t^2$ is bounded above by
$E_1+E_2+E_3$, \eqref{nse12a_gen} follows from the fact that we can choose
$K$ large enough to bound all three terms by some constant times
$\sigma^2\z(t)\abs{\log\eps}$. 
\end{proof}
\goodbreak

In the stable case, which is our major concern in this section, the
previous bound can be improved as follows.

\begin{cor}
\label{cor_nse}
Fix $t$ such that $t-t_0$ is at most of order $1$. In the stable case,
there exist constants $c_1>0$ and $M>0$ such that, if  
$\sigma \abs{\log\eps} \hat\z(t)^{3/2} \leqs c_1$, then
\begin{align}
\label{nse12a_spe}
\bigexpecin{t_0,0}{\abs{y_t}^{2k}} &\leqs k!\,M^k\sigma^{2k}\z(t)^k,
\qquad k\in\N.
\end{align}
\end{cor}
\begin{proof}
Let us again focus on $k=1$. Estimate~\eqref{nse12a_spe} is obtained
in the same way as~\eqref{nse12a_gen}, the only difference lying in a
more elaborate bound on $E_1$. We use the integral representation
\begin{equation}
\label{cor_nse:1}
y_t = \frac1\eps \int_{t_0}^t \e^{\alpha(t,s)/\eps} b(y_s,s)\6s +
\frac{\sigma}{\sqrt\eps} \int_{t_0}^t \e^{\alpha(t,s)/\eps} \6W_s
\end{equation}
of the SDE~\eqref{nse2} (for $y_{t_0}=0$), thereby obtaining
\begin{equation}
\label{cor_nse:3}
E_1\leqs 
2 \biggexpecin{t_0,0}
{\Bigpar{\frac1\eps \int_{t_0}^t \e^{\alpha(t,s)/\eps}
b(y_s,s)\6s}^2\,\bigindexfct{\sup_{t_0\leqs s\leqs
t}\frac{\abs{y_s}}{\sqrt{\z(s)}}\leqs\gamma}} 
+ 2 \frac{\sigma^2}\eps \int_{t_0}^t \e^{2\alpha(t,s)/\eps}\6s.
\end{equation}
The second term on the right-hand side is bounded above by
$2\sigma^2\z(t)$. The first one can be estimated by bounding
$b(y_s,s)$ uniformly in $s$, with the help of the estimate
$\abs{b(y,s)}\leqs M_0(y^2+\abs{y}^3)$, valid for all
$s$, c.f.~\eqref{nse3}. The remaining integral behaves like
$\eps \z(t)$. Thus we obtain
\begin{equation}
\label{cor_nse:4}
E_1 \leqs \text{\it const\/\,} \z(t)^2
\bigpar{\gamma^2\hat\z(t) + \gamma^3\hat\z(t)^{3/2}}^2 +
2\sigma^2\z(t),
\end{equation}
while $E$ can be estimated as before. Again choosing $\gamma =
K\sigma\sqrt{\abs{\log\eps}}$ for $K$ large, yields
Estimate~\eqref{nse12a_spe}. 
\end{proof}


\subsection{Noise-induced deviations from the deterministic area}
\label{ssec_areadev}

Let us now examine the behaviour of the surface delimited by the process
$y_t$. We want to control the process
\begin{equation}
\label{area1}
Y_t = - \int_{t_0}^t y_s \lambda'(s) \6s,
\end{equation}
which measures the deviation of the area $\cA(\eps,\sigma)$ enclosed
by $x_t$ from the one enclosed by $\xdet_t$.
Using the representation \eqref{cor_nse:1} of $y_t$, we obtain, by a
version of Fubini's theorem, that
\begin{equation}
\label{area3}
Y_t = \frac1\eps \int_{t_0}^t g(t,s) b(y_s,s)\6s + \frac\sigma{\sqrt\eps}
\int_{t_0}^t g(t,s)\6W_s,
\end{equation}
where 
\begin{equation}
\label{area4}
g(t,s) = - \int_s^t \e^{\alpha(u,s)/\eps}\lambda'(u)\6u.
\end{equation}
In particular, the term
\begin{equation}
\label{area5}
Y^0_t = \frac{\sigma}{\sqrt\eps} \int_{t_0}^t g(t,s)\6W_s
\end{equation}
is a Gaussian random variable with mean zero and variance
\begin{equation}
\label{area6}
\variance(Y^0_t) = \frac{\sigma^2}\eps \int_{t_0}^t g(t,s)^2 \6s.
\end{equation}
In the sequel, we will use the following abbreviations:
\begin{align}
\label{area6a}
\Gamma_i(t,t_0) & = \frac1\eps \int_{t_0}^t \abs{g(t,s)}\z(s)^i \6s, 
\qquad i\in\{1,\tfrac32\},\\
\label{area6b}
\Gamma(t,t_0) & = \frac1{\eps^2} \int_{t_0}^t \abs{g(t,s)}^2 \6s, \\
\label{area6d}
\Lambda(t,t_0) & = \int_{t_0}^t \abs{\lambda'(s)} \6s. 
\end{align}
In addition, we denote by $(2k-1)!!$ the product $\prod_{i=1}^{k} (2i-1)$.

\begin{prop}
\label{prop_area1}
Under the assumptions of Proposition~\ref{cor_nse_gen}, there is a constant
$M_1>0$ such that for all $k\geqs1$,
\begin{align}
\label{area7a1}
\bigexpecin{t_0,0}{(Y^0_t)^{2k}}
&= (2k-1)!!\, \bigpar{\sigma^2\eps \Gamma(t,t_0)}^k\\
\label{area7a2}
\bigexpecin{t_0,0}{\abs{Y^0_t}^{2k-1}}
&= k!\,2^{k-1} \sqrt{\frac2\pi} \bigpar{\sigma^2\eps
\Gamma(t,t_0)}^{(2k-1)/2}\\ 
\nonumber
\bigexpecin{t_0,0}{\bigabs{Y_t-Y^0_t}^k} 
&\leqs k!\, M_1^k \bigpar{\sigma^2\Gamma_1(t,t_0)}^k
\biggbrak{1 + \frac{\abs{\log\eps}^k}{k!}}\\
\label{area7b}
&\phantom{{}={}}\biggpar{1 + k!\,\sigma^{2k} \biggpar{\frac{\Gamma_{3/2}(t,t_0)}{\Gamma_1(t,t_0)}}^{2k}
\biggbrak{1 + \frac{\abs{\log\eps}^k}{k!}}}^{1/2}. 
\end{align}
\end{prop}
\begin{proof}
Since \eqref{area7a1} and \eqref{area7a2} are an immediate consequence
of the fact that $Y^0_t$ is Gaussian with variance~\eqref{area6}, we
only need to prove \eqref{area7b}. We restrict our attention to the
case $k=2$ as the case $k$ even follows by an obvious adaptation and
the case $k$ odd is obtained from the case $k$ even by an
application of Schwarz' inequality. First note that 
\begin{equation}
\label{p_area1:1new}
\bigexpecin{t_0,0}{(Y_t-Y^0_t)^2}
\leqs \frac{M_0^2}{\eps^2} \int_{t_0}^t \int_{t_0}^t \abs{g(t,u)}  \abs{g(t,v)}
 \bigexpecin{t_0,0}{(y_u^2+\abs{y_u}^3)(y_v^2+\abs{y_v}^3)}\6v\6u.
\end{equation}
Estimating the expectation of the product by H\"older's inequality and
Proposition~\ref{cor_nse_gen}, \eqref{area7b} follows.
\end{proof}
\goodbreak

\begin{remark}
\label{rem_area1}
In the stable case, under the assumptions of Corollary~\ref{cor_nse},
the bound~\eqref{area7b} simplifies to
\begin{equation}
\label{area7c}
\bigexpecin{t_0,0}{\bigabs{Y_t-Y^0_t}^k} 
\leqs k!\, M_1^k \bigpar{\sigma^2\Gamma_1(t,t_0)}^k
\bigpar{1 + k!\,\sigma^{2k}\hat\z(t)^k}^{1/2}. 
\end{equation}
\end{remark}

The following proposition gives bounds on the probability that the
deviation of the area from the corresponding area in the deterministic
case is large.

\begin{prop} \hfill
\label{prop_area2}
\begin{itemiz}
\item
Assume that there exists an $h_0(\eps,t)>0$ such that~\eqref{nse9} holds for 
all $h\leqs h_0(\eps,t)$. Then there exist constants $h_1, \kappa, C >
0$ such that for any $p\in(0,1)$,  
\begin{equation}
\label{area9}
\bigprobin{t_0,0}{\abs{Y_t}>H} \leqs 
\exp\biggset{-\frac{(1-p)^2}{2\Gamma(t,t_0)}\frac{H^2}{\sigma^2\eps}}
+C\biggpar{\frac{t-t_0}\eps + 1}
\exp\Bigset{-\frac{\kappa}{\sigma^2}\frac{p H}{\Gamma_1(t,t_0)}}, 
\end{equation}
whenever $p H\leqs h_1 \Gamma_1(t,t_0) (h_0(\eps,t)^2\wedge \hat\z(t)^{-1})$.
\item
In addition,
\begin{equation}
\label{area9a}
\bigprobin{t_0,0}{\abs{Y_t}>H} \leqs 
C\biggpar{\frac{t-t_0}\eps + 1} 
\exp\Bigset{-\frac{\kappa}{\sigma^2}
\Bigpar{\frac{H}{\Lambda(t,t_0)}}^{4}},
\end{equation}
whenever $H\geqs L_0\Lambda(t,t_0)$.
\end{itemiz}
\end{prop}
\begin{proof}
Consider first the case $pH$ small. By~\eqref{area3}, we have for
any $p\in(0,1)$  
\begin{equation}
\label{p_area2:1}
\bigprobin{t_0,0}{\abs{Y_t}>H} \leqs \bigprobin{t_0,0}{\abs{Y^0_t}>(1-p)H} + 
\Bigprobin{t_0,0}{\frac1\eps\int_{t_0}^t \abs{g(t,s)}\abs{b(y_s,s)}\6s > pH}.
\end{equation}
The first term on the right-hand side immediately yields the first term in
\eqref{area9} due to the Gaussian nature of $Y^0_t$. We denote $p H/M_0$
by $Q$. For $q\in(0,1)$, the second term can be bounded by 
\begin{equation}
\label{p_area2:a}
\biggprobin{t_0,0}{\sup_{t_0\leqs s\leqs t} \frac{\abs{y_s}}{\sqrt{\z(s)}} 
> \biggpar{\frac{qQ}{\Gamma_1(t,t_0)}}^{1/2}} 
+ \biggprobin{t_0,0}{\sup_{t_0\leqs s\leqs t}
\frac{\abs{y_s}}{\sqrt{\z(s)}} 
> \biggpar{\frac{(1-q)Q}{\Gamma_{3/2}(t,t_0)}}^{1/3}}. 
\end{equation}
We choose $q$ by
\begin{equation}
\label{p_area2:b}
\frac{q}{1-q} = \frac{1}{h_0(\eps,t)\wedge \hat\z(t)^{-1/2}}
\frac{\Gamma_1(t,t_0)}{\Gamma_{3/2}(t,t_0)} 
\end{equation}
and estimate both summands in~\eqref{p_area2:a} by~\eqref{nse9}. Note
that the first summand dominates the second one by our choice of $q$,
since we assumed $pH \leqs h_1 \Gamma_1(t,t_0) (h_0(\eps,t)^{2}\wedge
\hat\z(t)^{-1})$. Thus we obtain the bound~\eqref{area9}.

For $H$ large, we employ the trivial bound $\abs{Y_t} \leqs
(\sup_{t_0\leqs s\leqs t} \abs{y_s}) \Lambda(t,t_0)$ together
with Estimate~\eqref{nse10}, thereby obtaining~\eqref{area9a}.  
\end{proof}

Choosing $p = p(H,\eps)$ in~\eqref{area9} by
\begin{equation}
\label{area10p}
\frac{p}{(1-p)^2} = \frac {H}{2\kappa\eps}
\frac{\Gamma_1(t,t_0)}{\Gamma(t,t_0)}   
\end{equation}
yields the following corollary.

\begin{cor}
\label{cor_area}
There exist constants $h_1, \kappa, C > 0$ such that
\begin{equation}
\label{area10cor}
\bigprobin{t_0,0}{\abs{Y_t}>H} \leqs C\biggpar{\frac{t-t_0}\eps + 1} 
\exp\biggset{-\frac{(1-p)^2}{2\Gamma(t,t_0)}\frac{H^2}{\sigma^2\eps}}
\end{equation}
for all $H$ satisfying 
$(1-p)^2H^2 \leqs 2\kappa h_1 \eps \Gamma(t,t_0) 
(h_0(\eps,t)^2 \wedge \hat\z(t)^{-1})$. 
Here $p$ is defined by \eqref{area10p}. Furthermore, whenever
$H \leqs \text{const}\, \eps \Gamma(t,t_0)/\Gamma_1(t,t_0)$, then
$1-p$ is bounded away from zero.
\end{cor}

In order to complete the proof of Theorem~\ref{thm_smallamp}, we have to
control the function $g(t,s)$, defined in \eqref{area4}. This task is 
simplified by using the following lemma. 

\begin{lemma}
\label{lem_area}
Assume that $\abs{a'(u)}\leqs a_1(\eps)$ and 
$a(u)\leqs -c\sqrt{\eps a_1(\eps)}$ for all
$u$ in an interval $[s,t]$, where $c>0$ is independent of $a_0$
and $\eps$. Then there exists a constant $d>0$, independent of
$s,t,c$, such that  
\begin{equation}
\label{area11}
g(t,s) \asymp -\frac{\eps}{\abs{a(s)}} 
\biggbrak{\lambda'(s) +
\biggOrder{\frac{\eps}{\abs{a(s)}}}}
\end{equation}
whenever $t-s \geqs d\eps/\abs{a(s)}$. 
\end{lemma}
\goodbreak

\begin{proof}[\sc Proof of Theorem~\ref{thm_smallamp}] 
First note that Proposition~\ref{prop_nse1} establishes the
bound~\eqref{r9} on the probability of a sample path crossing the
potential barrier.
From \eqref{det10}, \eqref{nse6} and the preceding lemma, one easily obtains 
\begin{equation}
\label{area13new}
\Gamma_1(\tfrac12,-\tfrac12) = \bigOrder{\abs{\log(\abs{a_0}\vee\eps)}}
\quad\text{and}\quad
\Gamma(\tfrac12,-\tfrac12) \asymp 1 
\end{equation}
in the stable case, while $\Lambda(\tfrac12,-\tfrac12) \asymp 1$ is
trivial. Now the preceding results imply the stated estimates. 
\end{proof}


\section{The large noise regime}
\label{sec_lnr}

In this section, we consider those parameter regimes in which the noise
intensity $\sigma$ is large enough to allow for transitions from one
potential well to the other one, with a probability close to
$1$. Depending on the amplitude, there are
three cases to consider:
\begin{itemiz}
\item	$a_0\leqs-\eps$ and $\sigma\geqs\abs{a_0}^{3/4}$;
\item	$\abs{a_0}\leqs\eps$ and $\sigma\geqs\eps^{3/4}$;
\item	$a_0\geqs\eps$ and $\sigma\geqs\esa0^{1/2}$.
\end{itemiz}
Actually, we will need to assume that $\sigma \geqs K\abs{a_0}^{3/4}$,
$\sigma\geqs K\eps^{3/4}$ or $\sigma\geqs K\esa0^{1/2}$, respectively,
for some large constant $K$, but in order not to overburden notations,
we will assume that $K=1$ is a possible choice.

Recall that in the deterministic case, transitions are impossible if
$a_0\leqs\gamma_0\eps$, and occur only after time
$-\tc+c_0\eps^{2/3}\smash{a_0^{-1/6}}$ if $a_0\geqs\gamma_1\eps$. It turns
out that under the above conditions on $\sigma$, transitions are likely to
occur some time {\em before\/} the potential barrier reaches its minimal
height or even vanishes. For brevity, we shall only discuss the case
$\abs{a_0}\leqs\eps$ in detail, but the other cases can be investigated
similarly (since transitions occur early, they are not influenced by the
details of the bifurcation or avoided bifurcation). 


\subsection{The transition time}
\label{ssec_tt}

We assume $\abs{a_0}\leqs \eps$ unless stated otherwise. 
By symmetry, we may restrict our attention to a half-period, say
$t\in[-1/4,1/4]$. Let $x_t$ be the solution of the SDE~\eqref{SDE} starting
at time $t_0=-1/4$ in the upper well, i.\,e., near $x^\star_+(t_0)$. We
define the transition time as the stopping time 
\begin{equation}
\label{tt1}
\tau^0 = \inf \bigsetsuch{s>t_0}{x_s\leqs0} \in(t_0,\infty],
\end{equation}
when $x_s$ crosses the $t$-axis for the first time. The choice of $x_s=0$ is
purely for convenience, and the qualitative behaviour of $\tau^0$
remains the same if $0$ is replaced by any level between $-\xc-\delta$ and
$\xc+\delta$ as long as $\delta>0$ is chosen in such a way that
$f(x,t)\leqs 0$ holds for all $\abs{x}\leqs \xc+\delta$ and all $t$ in
question. The following result characterizes the distribution of $\tau^0$.

\begin{prop}
\label{prop_tt1}
There exist constants $C$, $c_1$, $c_2$, $\kappa > 0$ such that 
\begin{itemiz}
\item	for $t_0 < t \leqs -c_1\sigma^{2/3}$, 
\begin{equation}
\label{tt2}
\bigprobin{t_0,x_0}{\tau^0<t} \leqs C\biggpar{\frac{t-t_0}\eps + 1} 
\exp \Bigset{-\frac{\kappa}{\sigma^2\hat\z(t)^3}};
\end{equation}
\item	for $-c_1\sigma^{2/3} + c_2\eps \leqs t \leqs c_1\sigma^{2/3}$, 
\begin{equation}
\label{tt3}
\bigprobin{t_0,x_0}{\tau^0>t} \leqs \frac32  
\exp \Bigset{-\frac{\kappa}{\abs{\log\sigma}}\frac1\eps
\int_{-c_1\sigma^{2/3}}^t \abs{a(s)}\6s} + \e^{-\kappa/\sigma^2}.
\end{equation}
\end{itemiz}
Recall that $a(s)$ is the linearization of the drift term along
$\xdetof{+}_s$ as defined in~\eqref{det10}, and $\hat\z(t)$ is defined
in~\eqref{nse7}.
\end{prop}
\begin{proof}
The first part is a direct consequence of \eqref{nse9} with
$h=h_1\hat\z(t)^{-3/2}$, where $h_1$ is chosen sufficiently small that 
the relation $\abs{x_s-\xdet_s}\leqs h_1/\z(s)$ for all $s\in[t_0,t]$
implies that $x_s>0$ for these $s$. 

The second part is an application of Theorem~2.7 in \cite{BG2} (with
$h=\text{{\it const }}\sigma\abs{\log\eps}^{1/2}$). Note that the
theorem naturally extends to the case $0<a_0\leqs\eps$.
In fact, the integrand in \eqref{tt3} should be the curvature of the
potential at the deterministic solution tracking the saddle
$x^\star_0(t)$, but the curvature behaves like $\abs{a(t)}$, compare
\cite[Proposition~4.3]{BG2}. 
\end{proof}

The condition $\sigma\geqs\eps^{3/4}$ implies that
$\hat\z(t)^{-1}\asymp\abs{a(t)}\asymp\abs{t}$ for $t\leqs-c_1\sigma^{2/3}$,
and thus the exponent in \eqref{tt2} scales like $\abs{t}^3/\sigma^2$. The
integral in \eqref{tt3} behaves like $\sigma^{2/3}(t+c_1\sigma^{2/3})$. 

Proposition~\ref{prop_tt1} shows that the transition is likely to occur
close to time $t_1=-c_1\sigma^{2/3}$, which satisfies $\lambda(t_1) +
\lc \asymp \sigma^{4/3}$. Therefore, we should compare the area
$\cA(\eps,\sigma)$ to a reference area $\hat\cA$ given by 
\begin{equation}
\label{tt4}
\frac12 \hat\cA = \int_{t_0}^{t_1} \xdetof{+}_s (-\lambda'(s))\6s  
+ \int_{t_1}^{t_2} \xdetof{-}_s (-\lambda'(s))\6s,
\end{equation}
where $t_0=-1/4$, $t_1=-c_1\sigma^{2/3}$, $t_2=1/4$, and we denote by
$\xdetof{+}_s$ the deterministic solution starting in $x_0$, which tracks
$x^\star_+(s)$, and by $\xdetof{-}_s$ a deterministic solution tracking
$x^\star_-(s)$. It is easy to check that 
\begin{equation}
\label{tt5}
\hat \cA - \cA_0 \asymp -\sigma^{4/3}.
\end{equation}
(This relation does not depend on the initial conditions of $\xdetof{\pm}_s$,
as long as they are sufficiently close to $x^\star_+(t_0)$ or
$x^\star_-(t_1)$, respectively).  

\begin{remark}
\label{rem_tt}
Proposition~\ref{prop_tt1} also holds for $a_0<-\eps$, with the same
exponents. 

In the case $a_0>\eps$, $\hat\z(t)^{-1}$ behaves like $\abs{a(t)}$, given
by \eqref{det17}. The bound \eqref{tt2} holds for $t_0<t\leqs t_1 = -\tc
-c_1(\sigma^{2/3}\wedge \sigma^{4/3}\smash{a_0^{-1/2}})$, with the exponent
replaced by $-\kappa \smash{a_0^{3/4}}\abs{t+\tc}^{3/2}/\sigma^2$ if
$\sigma\leqs \smash{a_0^{3/4}}$ and by $-\kappa \abs{t+\tc}^{3}/\sigma^2$ if
$\sigma\geqs \smash{a_0^{3/4}}$. Note that in both cases,
$\lambda(t_1)+\lc\asymp \sigma^{4/3}$. 

The bound \eqref{tt3} holds for $t_1+c_2\eps \leqs t \leqs
-\tc+c_0\eps^{2/3}\smash{a_0^{-1/6}}$, with an exponent of the same order
as in the other cases, namely $\sigma^{2/3}(t-t_1)/(\eps\abs{\log\sigma})$.
The behaviour for larger $t$ will be discussed in
Proposition~\ref{prop_laa1} below.
\end{remark}


\subsection{Deviations from the reference area}
\label{ssec_dra}

Our aim is to characterize the deviations of the random variable
$\cA(\eps,\sigma)$ from its deterministic reference value $\hat\cA$
over one half-period. We focus again on the case $\abs{a_0}\leqs \eps$. 
With a slight abuse of notation, we can write this deviation as 
$\tfrac12(\cA(\eps,\sigma)-\hat\cA) = Y^+_{t_1} + Y^-_{t_2}$, where  
\begin{equation}
\label{dra1}
Y^+_{t_1} \defby \int_{t_0}^{t_1} (x_s - \xdetof{+}_s) (-\lambda'(s)) \6s, 
\qquad
Y^-_{t_2} \defby \int_{t_1}^{t_2} (x_s - \xdetof{-}_s) (-\lambda'(s)) \6s. 
\end{equation}
We will estimate separately the probability that each of these terms is
larger than $H$ or smaller than $-H$. To do so, we need a preparatory
result allowing to extend the estimate \eqref{nse9} to larger values of $h$.

\begin{prop}
\label{prop_dra1}
Define the stopping time 
\begin{equation}
\label{dra2}
\tau = \inf\bigsetsuch{t\in[t_0,t_1]}{x_t \leqs \xdetof{+}_t -
h_0\hat\z(t)^{-1}} \in[t_0,t_1]\cup\{\infty\},
\end{equation}
where the constant $h_0$ is taken from Proposition~\ref{prop_nse1}. Then 
\begin{equation}
\label{dra3}
\biggprobin{t_0,x_0}{\sup_{t_0\leqs s\leqs \tau}
\frac{\abs{x_s-\xdetof{+}_s}}{\sqrt{\z(s)}} > h} 
\leqs \frac C\eps \e^{-\kappa h^2/\sigma^2}
\end{equation}
for some $C$, $\kappa>0$, all $t\in[t_0,t_1]$ and\/ {\em all}
$h>0$. 
\end{prop}
\begin{proof}
The fact that the drift term $F$ has a negative second derivative with
respect to $x$ for all $x>0$ implies that $x_s$ is unlikely to exit the
strip of width $h\sqrt{\z(s)}$ through its upper boundary, as was
proved for negative $a_0$ in \cite[Proposition~4.5]{BG2}. We also know
by \eqref{nse9} that $x_s$ is unlikely to exit the strip through its
lower boundary if $\sigma\ll h\leqs h_0\hat\z(t)^{-3/2}$. The stopping
time $\tau$ has been defined in such a way that $x_s$ cannot leave a
strip of larger width before time $\tau$.  
\end{proof}

Note that by decreasing $h_0$ if necessary, we can arrange for
$\tau<\tau_0$. We are now able to estimate deviations of $Y^+_{t_1}$. 

\begin{prop}
\label{prop_dra2}
There exist constants $C$, $\kappa$, $h_1 > 0$ such that 
\begin{align}
\label{dra4}
\bigprobin{t_0,x_0}{Y^+_{t_1} < -H} &\leqs \frac C\eps \e^{-\kappa
H^{3/2}/\sigma^2} \\
\label{dra5}
\bigprobin{t_0,x_0}{Y^+_{t_1} > +H} &\leqs \e^{-\kappa H^2/(\sigma^2\eps)} +
\frac C\eps \e^{-\kappa H/(\sigma^2\abs{\log\sigma})}
\end{align}
for $0\leqs H\leqs h_1 \sigma^{2/3}\abs{\log\sigma}$. 
\end{prop}
\begin{proof}
We decompose 
\begin{align}
\label{p_dra2:1}
\bigprobin{t_0,x_0}{Y^+_{t_1} < -H} 
\leqs{} & \bigprobin{t_0,x_0}{Y^+_{\tau\wedge t_1} < -\tfrac12 H} \\
\nonumber
&+ \biggexpecin{t_0,x_0}{\indexfct{\tau<t_1}
\biggprobin{\tau,x_\tau}{\int_\tau^{t_1} (x_s-\xdetof{+}_s)
(-\lambda'(s))\6s < -\tfrac12 H}},
\end{align}
where $\tau$ is defined in~\eqref{dra2}. The first term on the right-hand
side can be estimated as in Proposition~\ref{prop_area2}, as there is no
need to distinguish positive and negative deviations for this term. However,
Proposition~\ref{prop_dra1} allows us to obtain bounds valid on a
larger domain of $H$. Note that Proposition~\ref{prop_area2} remains valid
when $Y_t$ is replaced by $Y_{\tau\wedge t}$. This is a consequence
of the $\Gamma_i(t,t_0)$ being monotone functions of $t\in[t_0,t_1]$ and a
slightly more elaborate estimate showing that $\sup_{t_0\leqs s \leqs
t} \abs{Y^0_s}$ obeys the same bound as was used for $\abs{Y^0_t}$
in~\eqref{p_area2:1}. Thus we obtain the estimate
\begin{equation}
\label{p_dra2:3}
\bigprobin{t_0,x_0}{\abs{Y^+_{\tau\wedge t_1}} > H} \leqs 
\exp\biggset{-\frac{(1-p)^2}{2\Gamma(t_1,t_0)}\frac{H^2}{\sigma^2\eps}} +
\frac C\eps\exp\Bigset{-\frac{\kappa}{\Gamma_1(t_1,t_0)}\frac{p H}{\sigma^2}}, 
\end{equation}
valid for $pH\leqs\text{{\it const }}\Gamma_1(t_1,t_0)/\hat\z(t_1)$.
An application of Lemma~\ref{lem_area} shows that $\Gamma(t_1,t_0)\asymp1$, 
$\Gamma_1(t_1,t_0)\asymp\abs{\log\sigma}$, and
we already know that $\smash{\hat\z(t_1)}\asymp\sigma^{-2/3}$. Choosing
$p\asymp1$ provides an estimate of the form \eqref{dra5}. 
The second term on the right-hand side of \eqref{p_dra2:1} can be estimated,
using the monotonicity of $\lambda$ in $[t_0,t_1]$, by the relation 
\begin{align}
\nonumber
& \Bigprobin{\tau,x_\tau}{\int_\tau^{t_1} (x_s-\xdetof{+}_s)(-\lambda'(s))\6s 
< -H}\\ 
\label{p_dra2:4}
&\qquad \leqs \Bigprobin{\tau,x_\tau}{\sup_{\tau\leqs s\leqs t_1}
\abs{x_s-\xdetof{+}_s}>L} 
+ \bigprobin{\tau,x_\tau}{\lambda(\tau)-\lambda(t_1) > H/L}. 
\end{align}
The first term on the right-hand side can be estimated by
Proposition~\ref{prop_nse2}, provided $L$ is larger than some constant of
order $1$. It decreases like $\e^{-\text{{\it const}}/\sigma^2}
\mskip-6mu/\eps$. Using the fact that 
$\lambda(\tau)-\lambda(t_1)\asymp\tau^2-t_1^2$ for 
$\tau$ not too close to $t_0$ (note that the contribution of $\tau$
close to $t_0$ is even smaller) and \eqref{tt2} of
Proposition~\ref{prop_tt1}, we obtain that  
\begin{equation}
\label{p_dra2:5}
\Bigexpecin{t_0,x_0}{\indexfct{\tau<t_1}
\bigprobin{\tau,x_\tau}{\lambda(\tau)-\lambda(t_1) > H/L}} 
\leqs \frac C\eps \e^{-\text{{\it const }}H^{3/2}/\sigma^2}.
\end{equation}
This last term is easily seen to dominate all others, so that
\eqref{dra4} is proved.  

To estimate deviations in the positive direction, we split terms as in
\eqref{p_dra2:1}. The first term can also be bounded by \eqref{p_dra2:3}.
Using the fact (compare \cite[Proposition~4.5]{BG2}) that
\begin{equation}
\label{p_dra2:6}
\biggprobin{\tau,x_\tau}{\sup_{\tau\leqs s\leqs t_1}
\frac{x_s-\xdetof{+}_s}{\sqrt{\z(s)}} > \sqrt{H/2}} 
\leqs \frac C\eps \e^{-\kappa H/\sigma^2}, 
\end{equation}
it only remains to estimate 
\begin{align}
\nonumber
&\biggprobin{\tau,x_\tau}{\int_\tau^{t_1} 
(x_s-\xdetof{+}_s)(-\lambda'(s))\6s > H/2, \,
\sup_{\tau\leqs s\leqs t_1} \frac{x_s-\xdetof{+}_s}{\sqrt{\z(s)}}
\leqs \sqrt{H/2}} \\
\nonumber
&\qquad \leqs \biggprobin{\tau,x_\tau}{\int_\tau^{t_1} \sqrt{\z(s)}
(-\lambda'(s)) \6s  > \sqrt{H/2}} \\
&\qquad\leqs \bigprobin{\tau,x_\tau}{\abs{\tau}^{3/2} > \sigma + \text{{\it
const }} \sqrt{H}}. 
\label{p_dra2:7}
\end{align}
The expectation of this term also decreases like $\e^{-\kappa
H/\sigma^2}\mskip-6mu/\eps$ by Proposition~\ref{prop_tt1}. Taking
\eqref{p_dra2:3} and \eqref{p_dra2:6} into account, we have proved
\eqref{dra5}.  
\end{proof}

The term $Y^-_{t_2}$ can be controlled in a similar way:

\begin{prop}
\label{prop_dra3}
There exist constants $C$, $\kappa$, $h_2$, $h_3 > 0$ such that 
for all $x_{t_1}\in[-L,L]$ ($L\asymp1$), and all $H\leqs h_3$,
\begin{align}
\label{dra6}
\bigprobin{t_1,x_{t_1}}{Y^-_{t_2} < -H} & \leqs 
\frac C\eps \e^{-\kappa H/\sigma^2} + 
\frac32 \e^{-\kappa\sigma^{4/3}/(\eps\abs{\log\sigma})} +  
\e^{-\kappa H^2/(\sigma^2\eps)} \\
\label{dra7}
\bigprobin{t_1,x_{t_1}}{Y^-_{t_2} > +H} 
&\leqs  \frac C\eps \e^{-\kappa H/\sigma^2} 
+ \frac32 \indicator{[0,h_2\sigma^{4/3})}(H) 
\e^{-\kappa H/(\eps\abs{\log\sigma})}. 
\end{align}
\end{prop}
\begin{proof}
The proof being similar to the one of the previous proposition, we only
outline the main steps. Introduce a stopping time
$\tauc=\inf\setsuch{s\in[t_1,t_2]}{x_s\leqs-\xc-\delta}$
for some small $\delta>0$, cf.~the comment on the definition of
$\tau^0$ in the beginning of the subsection. 
We first need to control the behaviour of 
\begin{equation}
\label{p_dra3:1}
Y^-_{\tauc\wedge t_2} = 
\int_{t_1}^{\tauc\wedge0} (x_s-\xdetof{-}_s) (-\lambda'(s))\6s + 
\int_{\tauc\wedge0}^{\tauc\wedge t_2} 
(x_s-\xdetof{-}_s) (-\lambda'(s))\6s. 
\end{equation}
Observe that since $x_s>\xdetof{-}_s$ for $s\leqs\tauc$, the first term on
the right-hand side is positive, while the second one is negative or zero.
First note that if $x_s$ is bounded above by $L\asymp1$, then
$Y^-_{\tauc\wedge t_2}$ cannot exceed a value of order
$\sigma^{4/3}$. Deviations of $Y^-_{\tauc\wedge t_2}$ in the positive
direction can be bounded using a decomposition similar to
\eqref{p_dra2:4} and applying \eqref{tt3} for $\tauc$ instead of
$\tau$. We find  
\begin{equation}
\label{p_dra3:2}
\bigprobin{t_1,x_{t_1}}{Y^-_{\tauc\wedge t_2}>H} \leqs 
\frac C\eps \e^{-\kappa /\sigma^2} 
+ \frac32 \indicator{[0,h_2\sigma^{4/3})}(H)  
\e^{-\kappa H/(\eps\abs{\log\sigma})},
\end{equation}
valid for $H\geqs\Order{\eps\sigma^{2/3}}$. 

Deviations of  $Y^-_{\tauc\wedge t_2}$ in the negative direction can
only be caused by the second term 
on the right-hand side of \eqref{p_dra3:1}. However, there is no small lower
bound for that term. The reason is that transitions to $\xdetof{-}_s$ are only
probable in the window $t\in[-c_1\sigma^{2/3},c_1\sigma^{2/3}]$. If this
opportunity is missed, which happens with a probability of order
$\e^{-\kappa\sigma^{4/3}/(\eps\abs{\log\sigma})}$, then $x_s$ keeps tracking
$\xdetof{+}_s$ and $Y^-_{\tauc\wedge t_2}$ may reach negative values of
order $1$. 

To complete the proof, we need to show that on $\{\tauc<t_2\}$
\begin{equation}
\label{p_dra3:3}
\biggprobin{\tauc,x_{\tauc}}{\biggabs{\int_{\tauc}^{t_2}
(x_s-\xdetof{-}_s)(-\lambda'(s))\6s} > H} \leqs 
\e^{-\kappa H^2/(\sigma^2\eps)} + \frac C\eps \e^{-\kappa H/\sigma^2}. 
\end{equation}
Let $\xdetof{\tauc}_s$ be the deterministic solution starting in
$x_{\tauc}$ at time $\tauc$. This solution is attractive, and thus
\eqref{p_dra3:3} holds with $\xdetof{\tauc}_s$ instead of $\xdetof{-}_s$ as a
consequence of Proposition~\ref{prop_area2}. But the distance between
$\xdetof{\tauc}_s$ and $\xdetof{-}_s$ decreases exponentially in
$(s-\tauc)/\eps$, which implies that the area between them is at most 
of order $\eps\sigma^{2/3}$. Thus \eqref{p_dra3:3} holds for
$H>\Order{\eps\sigma^{2/3}}$. But for smaller $H$, it is trivially
satisfied. 
\end{proof}
\goodbreak

We can summarize the properties obtained so far in the following way. 

\begin{prop}
\label{prop_dra4}
There exist constants $C$, $\kappa$, $h_1$, $h_2$, $h_3 > 0$ such that 
\begin{align}
\label{dra8}
\bigprobin{t_0,x_0}{\cA-\hat\cA < -H} &\leqs 
\frac C\eps \e^{-\kappa H^{3/2}/\sigma^2} + 
\frac32 \e^{-\kappa\sigma^{4/3}/(\eps\abs{\log\sigma})} \\
\bigprobin{t_0,x_0}{\cA-\hat\cA > +H} &\leqs 
\frac C\eps \e^{-\kappa H/(\sigma^2\abs{\log\sigma})} + 
\frac32 \indicator{[0,h_2\sigma^{4/3})}(H) \e^{-\kappa
H/(\eps\abs{\log\sigma})}
\label{dra9}
\end{align}
for $0\leqs H\leqs h_1\sigma^{2/3}\abs{\log\sigma}$. In addition, for all
$H\geqs h_3$ we have 
\begin{equation}
\label{dra10}
\bigprobin{t_0,x_0}{\abs{\cA}>H} \leqs \frac C\eps \e^{-\kappa
H^4/\sigma^2}.
\end{equation}
\end{prop}

As an immediate consequence, we obtain the following estimates on the
moments of the deviation of the area.

\begin{cor}
\label{cor_ev}
There exist positive constants $C$, $c_1$, $c_2$ such that 
\begin{align}
\label{dra12a}
\bigexpecin{t_0,x_0}{\cA(\eps,\sigma)-\hat\cA} 
&\geqs -C \sigma^{4/3} \abs{\log\eps}^{2/3} \\
\label{dra12b}
\bigexpecin{t_0,x_0}{\cA(\eps,\sigma)-\hat\cA} 
&\leqs C \bigbrak{(\eps \vee \sigma^2\abs{\log\eps})\abs{\log\sigma}}
\\
\label{dra12c}
\bigexpecin{t_0,x_0}{(\cA(\eps,\sigma)-\hat\cA)^2} 
&\leqs C \bigpar{\sigma^{4/3} \abs{\log\eps}^{2/3}}^2 
\end{align}
provided $c_1\eps<\sigma^{4/3}/\abs{\log\sigma}^2$ and
$\sigma^{2/3}\abs{\log\sigma} \leqs c_2/\abs{\log\eps}$. 
\end{cor}
\begin{proof}
By partial integration, 
\begin{equation}
\label{c_ev:1}
\bigexpecin{t_0,x_0}{\cA-\hat\cA} = 
\int_0^\infty \bigprobin{t_0,x_0}{\cA-\hat\cA>H}\6H - 
\int_0^\infty \bigprobin{t_0,x_0}{\cA-\hat\cA<-H}\6H. 
\end{equation}
The first integral be evaluated by splitting it at
$H=C(\eps \vee \sigma^2\abs{\log\eps})\abs{\log\sigma}$,
$h_2\sigma^{4/3}$ and $h_3$, and the second one 
at $C\sigma^{4/3}\abs{\log\eps}^{2/3}$ and $h_3$. Each time, the integral
over the first interval dominates. The estimate~\eqref{dra12c} is obtained
similarly. 
\end{proof}
\goodbreak


\section{The large amplitude case}
\label{sec_lac}

We consider finally the large amplitude case $a_0\geqs\gamma_1\eps$,
but with noise intensity $\sigma$ satisfying
$\sigma^2<\sqrt{a_0}\eps$. By symmetry, we may again concentrate on a
half-period $[-1/2,0]$. Let $x_t$ be the solution of the SDE
\eqref{SDE} starting at time $t_0=-1/2$ in the upper well, i.\,e.,
near $x^\star_+(t)$. Recall that the solution 
$\xdet_t$ of the deterministic equation \eqref{det1} with the same initial
condition tracks $x^\star_+(t)$ until time $-\tc$, and jumps to the other
potential well at $x^\star_-(t)$ after a delay of order
$\smash{\eps^{2/3}a_0^{-1/6}}$. 

We introduce a time $t^\star=-\tc+c_0\eps^{2/3}a_0^{-1/6}$ just before the
jump. Then we know that for $t_0\leqs t=-\tc+s\leqs t^\star$, 
\begin{equation}
\label{lac1}
\xdet_t-\xc \asymp \frac1{\z(t)} \asymp -a(t) 
\asymp \abs{s} \vee a_0^{1/4}\sqrt s \vee \esa0^{1/3},
\end{equation}
compare \eqref{det17} and \eqref{nse6}. The fact that $\xdet_t$ behaves in
this way follows from the fact that $x^\star_+(t)-\xc \asymp \sqrt{\mu(s)}$
dominates $x^\star(t)_+-\xdet_t$ for $t\leqs
-\tc-\smash{\eps^{2/3}a_0^{-1/6}}$, c.f.~\eqref{det13}. For larger $t$, we
know from \eqref{det15} and \eqref{det16} that $\xdet_t-\xc \asymp
\esa0^{1/3}$. 


\subsection{The transition time}
\label{ssec_lat}

Let us again start by investigating the distribution of the stopping time 
\begin{equation}
\label{lat1}
\tau^0 = \inf\bigsetsuch{s>t_0}{x_s=0}.
\end{equation}
The following result shows that $\tau^0$ is likely to be close to $t^\star$.

\begin{prop}
\label{prop_laa1}
There exist constants $C$, $c_1$, $\kappa>0$ such that 
\begin{itemiz}
\item	for $t_0\leqs t\leqs t^\star$, 
\begin{equation}
\label{lat2}
\bigprobin{t_0,x_0}{\tau^0<t} \leqs 
C\biggpar{\frac{t-t_0}\eps + 1} \exp\biggset{-\frac{\kappa}{\sigma^2
\hat\z(t)^3}}. 
\end{equation}

\item	for $t^\star \leqs t \leqs t^\star + c_1\sqrt{a_0}$, 
$0\leqs x_{t^\star}\leqs\rho$ and any $\rho>\xc$, 
\begin{equation}
\label{lat3}
\bigprobin{t^\star,x_{t^\star}}{\tau^0>t} \leqs 
3 \exp\biggset{-\frac{\eps\sqrt{a_0}}{\sigma^2}
\biggbrak{\frac{\kappa}{\abs{\log(\eps^2a_0)}\vee\log\rho} 
\frac{a_0^{1/6}(t-t^\star)}{\eps^{2/3}}-1}}.
\end{equation}
\end{itemiz}
\end{prop}
\begin{proof}
First note that~\eqref{lat2} is a direct consequence of
Proposition~\ref{prop_nse1} as we are in the stable case.

In order to prove~\eqref{lat3}, we consider again the stochastic
process $y_t=x_t-\xc$, satisfying the SDE 
\begin{equation}
\label{p_laa1:1}
\6y_t = \frac1\eps \bigbrak{\mu(t+\tc) - \sqrt3 y_t^2 - y_t^3} \6t 
+ \frac\sigma{\sqrt\eps} \6W_t, 
\qquad\qquad
y_{t^\star} = x_{t^\star}-\xc, 
\end{equation}
where \eqref{det13} implies $\mu(t+\tc)\leqs-\eps^{2/3}a_0^{1/3}$ for
$t^\star\leqs t\leqs t^\star+c_1\sqrt{a_0}$. Note furthermore that $- \sqrt3
y_t^2 - y_t^3\leqs -y_t^2$ for $t\leqs\tau^0$. By Gronwall's inequality, it
follows that $y_t\leqs z_t$ for $t\leqs \tau^0$, where $z_t$ is defined as
the solution of the time-homogeneous SDE
\begin{equation}
\label{p_laa1:2}
\6z_t = \frac1\eps \bigbrak{-\eps^{2/3}a_0^{1/3} - z_t^2} \6t 
+ \frac\sigma{\sqrt\eps} \6W_t, 
\qquad\qquad
z_{t^\star} = x_{t^\star}-\xc.
\end{equation}
For any $\delta_0>0$, we can write 
\begin{equation}
\label{p_laa1:3}
\bigprobin{t^\star,x_{t^\star}}{\tau^0>t} 
\leqs 
\Bigprobin{t^\star,x_{t^\star}}{\sup_{t^\star\leqs s\leqs t}z_s > \delta_0} 
+ 
\bigprobin{t^\star,x_{t^\star}}{-\xc\leqs z_s\leqs\delta_0 \;\forall
s\in[t^\star,t]}. 
\end{equation}
Since~\eqref{p_laa1:2} is an autonomous SDE, it is easy to see that
the first term on the right-hand side can be bounded by  
\begin{equation}
\label{p_laa1:4}
\Bigprobin{t^\star,x_{t^\star}}{\sup_{t^\star\leqs s\leqs t}z_s > \delta_0} 
\leqs C \Bigpar{\frac{t-t^\star}\eps+1} 
\e^{-\kappa\delta_0^3/\sigma^2},
\end{equation}
which can be made as small as we like by taking $\delta_0$ sufficiently
large. In order to estimate the second term, we introduce
$\Delta=c\eps^{2/3}a_0^{-1/6}$, where $c>1$ will be chosen later, and define
\begin{equation}
\label{p_laa1:5}
Q = \sup_{-\xc\leqs z_0\leqs\delta_0} \bigprobin{0,z_0}{-\xc\leqs
z_s\leqs\delta_0 \;\forall s\in[0,\Delta]}.
\end{equation}
Using time homogeneity and the Markov property, we can write
\begin{equation}
\label{p_laa1:6}
\bigprobin{t^\star,x_{t^\star}}{-\xc\leqs z_s\leqs\delta_0 \;\forall
s\in[t^\star,t]} 
\leqs Q^{(t-t^\star)/\Delta - 1}.
\end{equation}
The result is thus proved if we manage to bound $Q$ by a term exponentially
small in $\eps\sqrt{a_0}/\sigma^2$.

In order to estimate $Q$, it is convenient to introduce the process $\tilde
z_t = -z_t/(\eps^{1/3}a_0^{1/6})$, which obeys the SDE
\begin{equation}
\label{p_laa1:7}
\6\tilde z_t = \frac{a_0^{1/6}}{\eps^{2/3}} \bigbrak{1 + \tilde z_t^2} \6t 
- \frac\sigma{\eps^{5/6}a_0^{1/6}} \6W_t, 
\qquad\qquad
\tilde z_{t^\star} = -\frac{x_{t^\star}-\xc}{\eps^{1/3}a_0^{1/6}}.
\end{equation}
Let $\rho=\xc/(\eps^{1/3}a_0^{1/6})$ and
$\delta=\delta_0/(\eps^{1/3}a_0^{1/6})$. Using again Markov property
and time-homogeneity shows that $Q\leqs Q_1+Q_2+Q_3$,  where 
\begin{eqnarray}
\nonumber
Q_1 =& \mskip-8mu\displaystyle{\sup_{-\delta\leqs\tilde z_0\leqs-1}}\mskip-8mu
&\bigprobin{0,\tilde z_0}{\tilde z_s<-1\;\forall s\in[0,\Delta/3]} \\
\label{p_laa1:10}
Q_2 =&  \mskip-8mu\displaystyle{\sup_{-1\leqs\tilde z_0\leqs1}}\mskip-8mu
&\bigprobin{0,\tilde z_0}{\tilde z_s<1\;\forall s\in[0,\Delta/3]} \\
\nonumber
Q_3 =& \mskip-8mu \displaystyle{\sup_{1\leqs\tilde z_0\leqs\rho}}\mskip-8mu
&\bigprobin{0,\tilde z_0}{\tilde z_s<\rho\;\forall s\in[0,\Delta/3]}. 
\end{eqnarray}
Since $1+\tilde z^2\geqs 1\vee\abs{\tilde z}$, each term can be easily
estimated by comparison with an appropriate linear 
or $\tilde z$-independent equation. Consider for instance $Q_1$. We
know that $\tilde z_t$ lies above the solution $\tilde z^0_t$ of the linear SDE
\begin{equation}
\label{p_laa1:11}
\6\tilde z^0_t = -\frac{a_0^{1/6}}{\eps^{2/3}}\tilde z^0_t  \6t 
- \frac\sigma{\eps^{5/6}a_0^{1/6}} \6W_t, 
\qquad\qquad
\tilde z^0_0 = -\delta,
\end{equation}
the solution of which at time $\Delta/3$ is a  Gaussian random
variable with mean $-\delta\e^{-c/3}$ and variance
$(1-\e^{-2c/3})\sigma^2/(2\eps\sqrt{a_0}\mskip1.5mu)$. We can thus estimate
\begin{equation}
\label{p_laa1:12}
Q_1 
\leqs \bigprobin{0,-\delta}{\tilde z^0_{\Delta/3}<-1} 
\leqs \exp\Bigset{-\frac{\eps\sqrt{a_0}}{\sigma^2}
\frac{(1-\delta\e^{-c/3})^2}{1-\e^{-2c/3}}},
\end{equation}
provided $\delta\e^{-c/3}<1$, i.\,e., $c>3\log\delta$. Now, $Q_2$ and
$Q_3$ allow for similar bounds, and the result thus follows from
\eqref{p_laa1:4} and \eqref{p_laa1:6}, taking $\delta_0$ and $c$
sufficiently large.  
\end{proof}


\subsection{The case IIa}
\label{ssec_laa}

We now examine the process $Y_t$ defined in \eqref{area1}, which describes
deviations from the deterministic area. Using Lemma~\ref{lem_area} and
\eqref{lac1}, it is easy to check that 
\begin{equation}
\label{laa1}
\Gamma_1(t^\star,t_0) \asymp \abs{\log(\eps^{2/3}a_0^{-1/6})}
\quad\text{and}\quad
\Gamma(t^\star,t_0) \asymp 1,
\end{equation}
where $t_0=-1/2$. Applying Remark~\ref{rem_area1} and
Corollary~\ref{cor_area}, it is 
straightforward to check that the distribution of $Y_t$ is close to a
Gaussian with variance proportional to $\sigma^2\eps$. 

The situation changes, however, for $t>t^\star$, because the deterministic
solution crosses a zone of instability between $\xc$ and $-\xc$, compare
\eqref{nse3}. This instability causes a spreading of paths which we
will now analyse in more detail. Let us introduce times $t^\star_1$
and $t^\star_2$ such that 
\begin{equation}
\label{laa2}
\xdet_{t^\star_1} = \xc - c_1\esa0^{1/3}, 
\qquad\qquad
\xdet_{t^\star_2} = -\xc-\delta,
\end{equation}
where $c_1>0$ and $\delta<\xc$. Then $t^\star_1-t^\star$ and
$t^\star_2-t^\star$ are both of order $\eps^{2/3}a_0^{-1/6}$. 
We now proceed to determining the behaviour of $\z(t)$, defined in
\eqref{nse5}, which measures the spreading of paths around the deterministic
solution. 

\begin{prop}
\label{prop_laa2}
Let $z_t = \xc - \xdet_t$. 
Then there exist constants $C$, $K>0$ such that 
\begin{align}
\label{laa3a}
\z(t) &\asymp \esa0^{-1/3} &
&\text{for $t^\star \leqs t\leqs t^\star_1$} \\
\label{laa3b}
\z(t) &\asymp \esa0^{-5/3} z_t^4&
&\text{for $t^\star_1 \leqs t\leqs t^\star_2$} \\
\label{laa3c}
\z(t) &\leqs C \bigbrak{\esa0^{-5/3} \e^{-K(t-t^\star_2)/\eps}+1} &
&\text{for $t^\star_2 \leqs t\leqs 0.$} 
\end{align}
\end{prop}
\begin{proof}
\eqref{laa3a} follows by an elementary calculation. Next, consider the time
interval $t^\star_1 \leqs t\leqs t^\star_2$. The variable $z_t$ satisfies
the differential equation
\begin{equation}
\label{p_laa2:1}
\dtot{z_t}t = \frac1\eps \bigbrak{-\mu(t)+\sqrt3 z_t^2 - z_t^3}, 
\qquad\qquad
\mu(t) = \lc - A \cos(2\pi t).
\end{equation}
Note that $z_t$ is monotonously decreasing and $\mu(t)\leqs 0$ for the
times under consideration. The linearization of the drift term $F$ at $z_t$
is $a(t)=2\sqrt3 z_t - 3z_t^2$. It follows that for $t^\star_1 \leqs s\leqs
t\leqs t^\star_2$, 
\begin{equation}
\label{p_laa2:2}
\alpha(t,s) = \int_s^t a(u)\6u 
\leqs \eps \int_{z_s}^{z_t} \frac{2\sqrt3 z-3z^2}{\sqrt3 z^2-z^3} \6z 
= \eps \log \biggpar{\frac{z_t^2(\sqrt3-z_t)}{z_s^2(\sqrt3-z_s)}},
\end{equation}
and thus
\begin{equation}
\label{p_laa2:3}
\e^{\alpha(t,s)/\eps} \leqs \frac{z_t^2(\sqrt3-z_t)}{z_s^2(\sqrt3-z_s)} 
\asymp \frac{z_t^2}{z_s^2}.
\end{equation}
More careful estimates, based on the inequalities
\begin{equation}
\label{p_laa2:4}
\frac1\eps \bigbrak{\sqrt3 z_t^2 - z_t^3} \leqs \dtot zt 
\leqs \frac1\eps \bigbrak{ct + \sqrt3 z_t^2},
\end{equation}
show that $\e^{\alpha(t,s)/\eps}$ is also bounded below by a constant times
$(z_t/z_s)^2$.  Now $\z(t)$ can be computed in the same way, by performing
the change of variables $s\mapsto z_s$, yielding \eqref{laa3b}. Finally,
\eqref{laa3c} follows easily from the fact that we are again in the stable
case for $t\geqs t^\star_2$. 
\end{proof}
\goodbreak

We now return to the SDE \eqref{nse2} for $a(t)$ given by~\eqref{det17} and
$b(y,t)=F(y,t)-a(t)y$. Following the proof of \cite[Proposition 3.10]{BG2},
it is easy to establish~\eqref{nse9} for all $h\leqs h_0\esa0^{5/6}$ and
all $t$. The condition on $h$ stems from the fact that the linear term
$a(t)y_t$ should dominate the nonlinear term $b(y_t,t)$ for all
realizations $\w$ satisfying $\abs{y_t(\w)}\leqs h\sqrt{\z(t)}$ $\forall
t$. 

The condition on $h$ implies that we need to require
$\sigma\leqs\esa0^{5/6}$ for \eqref{nse9} to be of interest. Then the
maximal spreading of paths will typically be of order
$\sigma\sqrt{\z(t^\star_2)} \asymp \sigma\esa0^{-5/6}$. 

Since $\abs{a'(t)}$ is no longer bounded for $t^\star_1 \leqs t\leqs
t^\star_2$, we cannot apply Lemma~\ref{lem_area} to compute the integrals
\eqref{area6a} and \eqref{area6b}. However, using the same change of
variables as in the proof of Proposition~\ref{prop_laa2}, it is not
difficult to establish that 
\begin{equation}
\label{laa4}
\Gamma_1(0,t_0) \asymp \esa0^{-2/3}
\quad\text{and}\quad
\Gamma(0,t_0) \asymp \eps^{-2/3}a_0^{1/6},
\end{equation}
where again $t_0=-1/2$.

The following proposition now follows immediately from
Corollary~\ref{cor_area}. 

\begin{prop}
\label{prop_laa3}
Assume that $\sigma\leqs\esa0^{5/6}$. There exists a constant $h_1$ such
that for $H\leqs h_1\eps\sqrt{a_0}$, 
\begin{equation}
\label{laa5}
\bigprobin{t_0,x_0}{\abs{Y_0}>H} \leqs \frac C\eps
\exp\Bigset{-\kappa\frac{H^2}{\sigma^2\esa0^{1/3}}}.
\end{equation}
\end{prop}

Finally, Proposition~\ref{prop_area1} can also be applied to show that for
$\sigma\abs{\log\eps} \leqs \text{{\it const }} \esa0^{5/6}$, 
\begin{align}
\label{laa6a}
\bigabs{\bigexpecin{t_0,x_0}{Y_0}} &=
\BigOrder{\frac{\sigma^2\abs{\log\eps}}{\esa0^{2/3}}} \\
\label{laa6b}
\variance{(Y_0)} &\asymp \sigma^2\esa0^{1/3}.
\end{align}


\subsection{The case IIb}
\label{ssec_lab}

For $\sigma$ larger than $\esa0^{5/6}$, the strong dispersion of
trajectories near time $t^\star_2$ prevents us from applying methods of
Section~\ref{sec_near}. However, methods similar to those
of Section~\ref{sec_lnr} can be applied to obtain some information. 

The following result allows to estimate deviations of $Y_{t^\star}$ in a
larger domain than Proposition~\ref{prop_area2}.

\begin{prop}
\label{prop_lab1}
There exist constants $C$, $\kappa$, $h_1 > 0$ such that 
\begin{align}
\label{lab1}
\bigprobin{t_0,x_0}{Y_{t^\star} < -H} &\leqs \frac C\eps \e^{-\kappa
(H^{3/2}\vee\eps\sqrt{a_0}\mskip1.5mu)/\sigma^2} \\
\label{lab2}
\bigprobin{t_0,x_0}{Y_{t^\star} > +H} &\leqs \e^{-\kappa H^2/(\sigma^2\eps)} +
\frac C\eps \e^{-\kappa H/(\sigma^2\Gamma_1(t^\star,t_0))}
\end{align}
for $0\leqs H\leqs h_1 \esa0^{1/3} \Gamma_1(t^\star,t_0)$, where 
$\Gamma_1(t^\star,t_0) \asymp \abs{\log(\eps^{2/3}a_0^{-1/6})}$. 
\end{prop}
\begin{proof}
The proof is almost the same as the proof of Proposition~\ref{prop_dra2},
the only difference lying in the different behaviour of $\z(t)$, given in
\eqref{lac1}, which requires to distinguish between
$\tau+\tc\leqs-\sqrt{a_0}$,
$-\sqrt{a_0}\leqs\tau+\tc\leqs-\eps^{2/3}\smash{a_0^{-1/6}}$, and the
remaining $\tau$ up to $t^\star$. 
\end{proof}
\goodbreak

Proceeding as in the proof of Proposition~\ref{prop_dra3}, but using
Proposition~\ref{prop_laa1} for the transition time, we obtain

\begin{prop}
\label{prop_lab2}
There exist constants $C$, $\kappa$, $h_2$, $h_3 > 0$ such that 
for all $x_{t_1}\in[-L,L]$ ($L\asymp1$), and
$h_2\esa0^{2/3}\abs{\log\esa0}\leqs H\leqs h_3$,
\begin{align}
\label{lab3}
\bigprobin{t^\star,x_{t^\star}}{Y_0 < -H} &\leqs 
\frac C\eps \e^{-\kappa H/\sigma^2} + 
\e^{-\kappa H^2/(\sigma^2\eps)} \\
\label{lab4}
\bigprobin{t^\star,x_{t^\star}}{Y_0 > +H} &\leqs 
\frac C\eps \e^{-\kappa \esa0^{1/3}H/(\sigma^2\abs{\log\esa0})}. 
\end{align}
\end{prop}

\begin{cor}
\label{cor_lab}
For $h_2\esa0^{2/3}\abs{\log\esa0}\leqs
h_1\esa0^{1/3}\Gamma_1(t^\star,t_0)$, 
\begin{align}
\label{lab5}
\bigprobin{t_0,x_0}{Y_0 < -H} &\leqs \frac C\eps \e^{-\kappa
H^{3/2}/\sigma^2} \\
\label{lab6}
\bigprobin{t_0,x_0}{Y_0 > +H} &\leqs 
\frac C\eps \e^{-\kappa \esa0^{1/3}H/(\sigma^2\abs{\log\esa0})}. 
\end{align}
\end{cor}

The required lower bound on $H$ only allows us to conclude that expectation
and standard deviation of $Y_0$ are smaller than a constant times
$\esa0^{2/3}\abs{\log\esa0}$, although the above estimates are already very
small for $H = h_2\esa0^{2/3}\abs{\log\esa0}$. 


\goodbreak

\bigskip\bigskip\noindent
{\small 
Nils Berglund \\ 
{\sc Department of Mathematics, ETH Z\"urich} \\ 
ETH Zentrum, 8092~Z\"urich, Switzerland \\
{\it E-mail address: }{\tt berglund@math.ethz.ch}

\bigskip\noindent
Barbara Gentz \\ 
{\sc Weierstra\ss\ Institute for Applied Analysis and Stochastics} \\
Mohrenstra{\ss}e~39, 10117~Berlin, Germany \\
{\it E-mail address: }{\tt gentz@wias-berlin.de}
}


\end{document}